\documentclass[a4paper,11pt,book]{amsart}

\usepackage[francais]{babel}
\usepackage{amssymb} 
\usepackage[latin1]{inputenc}
\usepackage[matrix,arrow]{xy}
\usepackage[mathscr]{eucal}
\usepackage{epic,eepic}
\usepackage{bbm}

\newtheorem{prop}[subsubsection]{Proposition}

\newtheorem{proposition}[subsubsection]{Proposition}
\newtheorem{theoreme}[subsubsection]{Théorème}

\newtheorem{lemme}[subsubsection]{Lemme}

\newtheorem{surtheoreme}{Théorème}
\newtheorem{surcorollaire}[surtheoreme]{Corollaire}

\theoremstyle{definition}
\newtheorem{definition}{Définition}

\newtheorem{remarque}{Remarque\!\!}

\newenvironment{pf}
        {\medskip\noindent {\it Démonstration --- \ }}
        {\hfill\nobreak $\Box$ \par\bigbreak}

\newcommand{\anti}{{\overline\ast}}

\newcommand{\Hom}{\text{Hom}}

\newcommand{\p}{{\mathfrak P}}

\newcommand{\C}{{ \mathbb C  }}
\newcommand{\R}{{ \mathbb R  }}
\newcommand{\Q}{{ \mathbb Q } }
\newcommand{\Z}{{ \mathbb Z  }}

\newcommand{\Ker}{{\mathrm{Ker}\,}}

\renewcommand{\ker}{{\mathrm{Ker}\,}}
\newcommand{\Gal}{{\mathrm{Gal}\,}}

\newcommand{\X}{{\mathcal X}}
\newcommand{\A}{\mathbb A}

\newcommand{\anneau}{{ \mathcal O}}

\newcommand{\Id}{{\mathrm{Id}}}

\newcommand{\U}{{\mathrm{U}}}

\newcommand{\Gl}{{\mathrm {GL}}}

\newcommand{\GL}{{\mathrm {GL}}}

\newcommand{\Sl}{{\mathrm {SL}}}

\newcommand{\sat}{{ \mathrm{sat}}}
\newcommand{\lng}{{ \mathrm{long}}}
\newcommand{\const}{{\mathop{\mathrm{const}}}}

\newcommand{\Zc}{{\mathcal Z}}

\newcommand{\spec}{{\mathrm{Spec\,}}}

\newcommand{\CC}{{\mathcal{C}}}

\newcommand{\St}{\mathrm{St}}
\newcommand{\st}{\mathrm{St}}
 
\newcommand{\Hecke}{{\mathcal{H}}}

\newcommand{\Sc}{{\mathcal{S}}}

\newcommand{\BB}{{\mathcal{B}}}

\newcommand{\pig}{{\varpi}}

\newcommand{\diag}{{\mathrm{diag}}}
\newcommand{\Frob}{{\mathrm{Frob\,}}}

\newcommand{\conv}{\ast}
\newcommand{\muz}{{\mu_0}}

\newcommand{\sd}{\rtimes}
\newcommand{\val}{{\mathrm{val}}}
\newcommand{\mat}{{\mathrm{mat}}}
\newcommand{\CCC}{{\mathcal{CC}}}
\newcommand{\vz}{{v_0}}
\newcommand{\ds}{\displaystyle}

\begin{document}

\baselineskip 14pt

\title{Augmentation du niveau pour $\U(3)$}
\author{Joël Bellaïche et Philippe Graftieaux}
\maketitle
\setcounter{tocdepth}{1} 
{\bf Résumé} : Nous démontrons pour la forme compacte à l'infini du
 groupe unitaire à trois variables attaché à une extension CM un
 résultat d'augmentation du niveau analogue à celui obtenu par Taylor 
(\cite{taylor}) 
dans le cas de $\Gl_2$. Nous donnons une application aux
 représentations automorphes non tempérées.

{\bf Abstract}: We prove, for the unitary group in three variables
attached to a CM extension which is compact at infinity, a
level-raising theorem analogous to the one of Taylor (\cite{taylor}) 
in the case of 
$\Gl_2$. We give an application to non tempered automorphic forms.
 

\tableofcontents
\bibliographystyle{style} 

\section{Introduction}

Soit $E/F$ une extension CM, et $G=U(3)$ le groupe unitaire à trois
variables sur $F$ attaché à $E$ qui est compact à toutes les
places à l'infini (cf.~\ref{groupeunitaire}). Pour toute
place finie $v$ de $F$, on note $F_v$ le complété de $F$ en $v$
et $G_v=G(F_v)$.   

Soient $K=\prod_{v \text{ place finie de $F$}} K_v$ un sous-groupe 
compact ouvert de $G(\A_{F,f})$
({\it le niveau}), $\Sigma$ l'ensemble fini
des places $v$ où $K_v$ n'est pas un sous-groupe compact
maximal hyperspécial, $J$ une représentation complexe lisse
irréductible de $\prod_{v \in \Sigma} K_v$ ({\it le type}), 
et $\rho$ une représentation complexe continue irréductible
de $G(\A_{F,\infty})$ ({\it le poids}). Pour tout ensemble fini $\Sigma'$
de places de $F$ contenant $\Sigma$, 
soit $\Hecke^{\Sigma'}$ l'algèbre de Hecke $
 \Hecke(\textstyle{\prod_{v \not \in \Sigma'}' G_v
,\prod_{v \not \in \Sigma'} K_v}).$

Soient $v_0 \not \in \Sigma$ une place finie de $F$  {\bf inerte} 
dans $E$, $q$ le cardinal du corps résiduel de $F_{v_0}$, 
$B_{v_0}$ un sous-groupe d'Iwahori de $K_{v_0}$, 
$B=B_{v_0} \prod_{v \neq v_0} K_v$ et $T_{v_0}\in \Hecke(G_{v_0},K_{v_0})$ 
l'opérateur de Hecke standard en $v_0$.

Nous notons $S_{K,J,\rho,\C}$ (resp. $S_{B,J,\rho,\C}$) 
l'espace vectoriel complexe des formes automorphes pour
$G$ de niveau $K$ (resp. $B$), type $J$ et poids $\rho$, et 
$O_{B,J,\rho,\C}$ (resp. $N_{B,J,\rho,\C}$)
le sous-espace de $S_{B,J,\rho,\C}$ des formes 
anciennes (resp. nouvelles) en $v_0$ 
(cf.~\ref{parformescomplexes} et~\ref{oldetnew}). 
Ces espaces sont munis d'une action naturelle de $\Hecke^\Sigma$. 
En particulier, on dispose d'une décomposition en sous-espaces propres
généralisés (cf.~\ref{sepgeneralises})
\begin{eqnarray}\label{decompintro}
S_{K,J,\rho,\C} &=& \oplus_\eta S_{K,J,\rho,\C}(\eta),
\end{eqnarray}
où $\eta$ décrit un ensemble fini de caractères complexes de
$\Hecke^\Sigma$, ainsi que d'une décomposition
analogue de $O_{B,J,\rho,\C}$.

Il existe un corps de nombres $L \subset \C$ et un
ensemble fini $S$ de place de $L$ ($L$ et $S$ ne dépendant
que de $(K,J,\rho)$ et pas de $v_0$) tels que pour toute place finie
$\mu$ de $L$ hors de $S$, l'espace $S_{B,J,\rho,\C}$
(resp. $O_{B,J,\rho,\C}$, $N_{B,J,\rho,\C}$) admette un
modèle $S_{B,J,\rho,\anneau_\mu}$
(resp. $O_{B,J,\rho,\anneau_\mu}$, $N_{B,J,\rho,\anneau_\mu}$)
sur $\anneau_\mu$, stable sous l'action de $\Hecke^{\Sigma_\mu}$, où
$\Sigma_\mu$ désigne la réunion de $\Sigma$ et des places finies $v$ de $F$ 
de caractéristique résiduelle distincte de celle de $\mu$, 
et $\anneau_\mu$ l'anneau des entiers 
du complété $L_\mu$ de $L$ en $\mu$. 
Les caractères $\eta$ intervenant en~(\ref{decompintro}) sont
alors à valeurs dans $\anneau_\mu$ 
et on a une décomposition 
$O_{B,J,\rho,\anneau_\mu}
=\oplus_{\eta} O_{B,J,\rho,\anneau_\mu}(\eta)$.

\begin{surtheoreme}\label{augniv} 
Soient $\psi$ un caractère de $\Hecke^\Sigma$ tel que 
$S_{K,J,\rho,\C}(\psi)$ est non nul et $\mu$ une place de $L$.
Supposons que $\mu$ est première à $v_0$ et que 
$\lambda:=\psi(T_{v_0}) \in \anneau_\mu$ vérifie $\lambda \neq q(q^3+1)$. 
Notons  
$c$ le plus petit entier vérifiant 
$$c \geq \val_\mu(\lambda-q(q^3+1))/2.$$ 
Il existe alors une congruence\footnote{i.e. il existe $f \in
O_{B,J,\rho,\anneau_\mu}(\psi)$, $g \in N_{B,J,\rho,\anneau_\mu}$, avec
$f-g \in \mu^c S_{B,J,\rho,\anneau_\mu}$, 
$f \not \in \mu S_{B,J,\rho,\anneau_\mu}$.}
modulo $\mu^c$ entre  
$O_{B,J,\rho,\anneau_\mu}(\psi)$ et $N_{B,J,\rho,\anneau_\mu}$
dans $S_{B,J,\rho,\anneau_\mu}$. 
\end{surtheoreme}
Nous renvoyons le lecteur au théorème~\ref{theoaugniv} et 
au paragraphe~\ref{remarques} pour des compléments à cet énoncé.

Le cas $c=1$ du théorème~\ref{augniv}, 
combiné avec le fameux lemme de Deligne-Serre,
donne en particulier le résultat suivant.
\begin{surcorollaire}\label{coraugniv}
Soient $\pi$ une représentation automorphe pour $G$, de dimension infinie 
et non ramifiée en $v_0$, et 
$\lambda$ la valeur propre de $T_{v_0}$ sur $\pi_{v_0}$. 
Il existe un corps de nombres
$L$ et un ensemble finie de place $S$ de $L$ 
(ne dépendant que de $\pi$) 
tel que si $\mu \not \in S$ est une place finie de $L$ vérifiant  
$$\lambda \equiv q(q^3+1) \pmod{\mu},$$
alors il existe une représentation automorphe $\pi'$, de même poids 
et de même niveau hors $v_0$ que $\pi$, vérifiant
$\pi_{v_0}^{B_{v_0}} \neq 0$ mais $\pi_{v_0}^{K_{v_0}}=0$ et qui est 
congrue\footnote{i.e. les polynômes caractéristiques des 
matrices de Hecke de $\pi$ et de $\pi'$
en toute place non ramifiée sont à coefficients dans $L$
et congrus modulo $\mu$} 
à $\pi$ modulo $\mu$. De plus, si $\pi$ contient un certain
type hors $v_0$, on peut aussi supposer qu'il en va de même pour $\pi'$.
\end{surcorollaire}

Ce corollaire généralise \cite[théorème 2.4]{clozel}, qui traite le
cas\footnote{À vrai dire, Clozel énonce son théorème pour la forme
quasi-déployée du groupe unitaire $G$. Mais c'est bien le théorème analogue
pour le groupe compact qui est prouvé, et Clozel en déduit son énoncé 
pour la forme quasi-déployée par un argument de
transfert entre les formes intérieures.}
où le poids de $\pi$ et le type considéré 
sont triviaux, et qui surtout
suppose que $\mu$ est une place {\it banale} pour $G_{v_0}$, i.e. 
$(q-1)(q^3+1)\not \equiv 0 \pmod{\mu}$. 
Dans~\cite[Théorème VII.1.4.6]{Joelthese}, 
l'hypothèse de banalité de $\mu$
avait été affaiblie en une hypothèse de {\it normalité} 
$(q^3+1) \neq 0 \pmod{\mu}$. 

La méthode que nous
 utilisons pour prouver le théorème n'est pas la même que celle 
de~\cite{clozel} (ou 
de \cite{Joelthese}) basée sur les propriétés du {\it
 module universel} et sur un argument de densité, inspiré de Serre, jouant
le rôle du {\it lemme d'Ihara}. Il semble aux auteurs que cette méthode 
ne peut permettre de lever l'hypothèse de {\it normalité}  
de $\mu$, ni
d'obtenir des congruences modulo $\mu^c$ avec $c >1$. Notre méthode 
se rapproche au contraire de celle de~\cite{taylor}, avec
 quelques différences importantes. La principale est que le lemme
 d'Ihara (\cite[lemme 4]{taylor}) est {\it faux} pour $G=\U(3)$ 
en caractéristique non  normale, même pour des formes de poids
 trivial. Nous devons en quelque sorte majorer son  défaut 
(lemme~\ref{unautrelemme}), à l'aide d'une étude combinatoire d'opérateurs
entres espaces de fonctions sur l'arbre bihomogène attaché à $G_{v_0}$
(lemme~\ref{relations}). 
Le fait de traiter des poids quelconques, 
qui ne sont pas nécessairement 
auto-duaux au caractère central près comme dans la situation 
de~\cite{taylor}, nécessite aussi quelques arguments supplémentaires. 

Notons enfin que le théorème~\ref{augniv} 
et son corollaire restent valables, avec
la même preuve, pour n'importe quel groupe algébrique connexe $G$ sur
$F$, compact aux places à l'infini et de rang un en la place $v_0$ de
$F$. Il suffit de remplacer dans l'énoncé du théorème et de son
corollaire $q(q^3+1)$  par $q^{d'}(q^d+1)$, où $q^{d'}+1$ et $q^d+1$
sont les valences de l'arbre de Bruhat-Tits\footnote{Rappelons
(\cite{tits}) que cet arbre est soit homogène,
soit bihomogène. On ramène le premier cas au second en introduisant un
sommet au milieu de chaque arête, ce qui revient à poser $d'=0$.} 
de $G_{v_0}$, avec $d'\leq d$. 

Revenons au cas $G=U(3)$. À l'aide de la représentation galoisienne
attachée par Blasius et Rogawski à une représentation automorphe $\pi$ de
$G$, on peut montrer qu'il existe une infinité de places $v_0$ où
l'on peut augmenter le niveau de $\pi$. Plus précisément :
\begin{surtheoreme} \label{thexistence}
Gardons les notations du théorème~\ref{augniv}.
Pour tout entier $n$, il existe un ensemble de densité 
non nulle de places $v$ de
$F$ inertes dans $E$ 
telles que $\psi(T_v) \equiv q(q^3+1) \pmod{\mu^n}$ et $q+1 \equiv 0
\pmod{\mu^n}$. 
\end{surtheoreme}
En combinant les théorèmes~\ref{augniv} et~\ref{thexistence}, on obtient 
que les réductions modulo $\mu^n$ du caractère $\psi$
apparaissent dans des espaces de formes nouvelles
pour presque tout $\mu$ et pour $n$ arbitrairement grand.
À cause de l'hypothèse de banalité, Clozel ne pouvait démontrer un tel
résultat (avec $n=1$) que sous une hypothèse de surjectivité de la
représentation galoisienne attachée à $\psi$ modulo $\mu$.   

Comme application arithmétique du théorème~\ref{augniv}, on obtient
que toute représentation endoscopique 
non tempérée de $G$ est congrue en presque
toute place à une représentation tempérée (voir le
théorème~\ref{thtemp}). 
Les congruences entre formes 
endoscopiques non tempérées et formes
tempérées ont une signification arithmétique importante, en ce
qu'elles traduisent et permettent de montrer des cas des conjectures de
Bloch-Kato.
Dans ce but, une version plus faible du théorème précédent, 
valable seulement pour
un ensemble de densité non nulle de places, était obtenue
par une augmentation
du niveau en une place décomposée dans~\cite[chapitre VIII]{Joelthese}, par
l'utilisation de familles $l$-adiques de représentations automorphes 
dans~\cite{joeletgaetan}. 
Obtenir ce résultat pour presque toute place était une des
motivations initiales de ce travail.

Dans un article en préparation, nous montrons comment on peut utiliser
le théorème~\ref{augniv} 
pour montrer en toute place la compatibilité (à semi-simplification
près) de la construction de Blasius-Rogawski (\cite{br})
d'une représentation galoisienne
attachée à une repr\'esentation automorphe pour $G$ avec la correspondance de
Langlands locale (compatibilité qui n'est connue jusqu'à présent 
que pour les places non ramifiées,
ou bien en toute place si l'on suppose que la représentation
automorphe a son changement de base à $E$ de carré intégrable en au
moins une place finie, d'après les travaux de Harris et Taylor \cite{ht})

{\bf Remerciements :} Les auteurs remercient chaleureusement 
Laurent Clozel et Gaëtan Chenevier 
pour de nombreuses et éclairantes conversations.  
\section{Notations}

Sauf mention explicite du contraire, tous les anneaux ou algèbres
considérés sont commutatifs et unitaires.

\subsection{Corps}
\label{corps}
Dans tout l'article,  $E/F$ désigne une extension CM, 
$\anneau_E$ (resp. $\anneau_F$) l'anneau des entiers de $E$ (resp. $F$), 
et $c$ l'élément non trivial de $\Gal(E/F)$. On fixe une clôture algébrique
$\bar E$ de $E$. Pour toute représentation $\rho$ de $\Gal(\bar E/E)$, on 
note $\rho^c$ la représentation $g \mapsto \rho(\gamma g \gamma^{-1})$, 
où $\gamma$ est un  relevé de $c$ dans $\Gal(\bar E/F)$; 
la représentation $\rho^c$ ne dépend du choix de 
$\gamma$ qu'à isomorphisme près.  

La lettre $v$ (resp. $w$), éventuellement munie d'indices, désigne une 
place finie de $F$ (resp. de $E$); on note $F_v$ (resp. $E_w$)
le complété de $F$ en $v$ 
(resp. de $E$ en $w$), d'anneau d'entiers $\anneau_v$ (resp. $\anneau_w$). 
La lettre $\sigma$ désigne une place archimédienne de $F$, 
ce que l'on note $\sigma|\infty$.

On note $\A_F$ (resp. $\A_{F,f}$, $\A_{F,\infty}$) l'anneau des adèles
de $F$
(resp. les sous-anneaux des adèles triviaux à l'infini, aux places
finies).

La lettre $L$ désigne un sous-corps de $\C$, 
dont les places sont désignées par la lettre $\mu$. 
On note $L_0$ le corps $L\cap \R$ et $L_\mu$ le complété de $L$ en $\mu$, 
d'anneau d'entiers $\anneau_\mu$.

\subsection{Groupe unitaire}\label{groupeunitaire}

On note $G$ l'unique groupe 
unitaire à trois variables sur $F$, 
compact à l'infini associé à l'extension $E/F$,
qu'on peut définir ainsi : pour toute $F$-algèbre $R$, on pose 
$G(R)=\{g \in \Gl_n(E \otimes_{F} R),\ c(g) g=1\}$.
On note encore $G$ le modèle de $G$ sur $\anneau_F$ obtenu en posant
$G(R)=\{g \in \Gl_n(\anneau_E \otimes_{\anneau_F} R),\ c(g) g=1\}$
pour toute $\anneau_F$-algèbre $R$.
  
Pour toute place finie $v$ de $F$, 
on note $G_v$ le groupe $G(F_v)$ et $Z_v$ son centre.

\subsection{Applications linéaires localement finies} 
Soit $R$ un anneau. 
Pour tout ensemble $X$, 
on note $\CC(X,R)$ (resp. $\CCC(X,R)$)
le $R$-module des fonctions sur $X$ à valeurs dans $R$ 
(resp. nulles en dehors d'un sous-ensemble fini de $X$).  
Nous munissons $R$ de la topologie discrète et 
$\CC(X,R)$ de la topologie produit,
de sorte que $\CCC(X,R)$ est dense dans $\CC(X,R)$. 
Pour tout $x\in X$, on note $\delta_x$
la fonction caractéristique de $\{x\} \subset X$. 

Si $X$ et $Y$ sont deux ensembles, on 
dit qu'une application linéaire $U\colon \CC(X,R) \to \CC(Y,R)$ 
est {\it localement finie} si elle est continue et prolonge
une application linéaire $\CCC(X,R) \to \CCC(Y,R)$. 
Ceci équivaut à dire 
que la famille $U(\delta_x)_{x\in X}$ est à valeurs dans  $\CCC(Y,R)$ 
et détermine $U$. 

Si $U\colon \CC(X,R) \to \CC(Y,R)$ 
est une application linéaire localement finie, 
on définit sa transposée $U^\ast\colon \CC(Y,k)\to \CC(X,k)$
comme l'application linéaire localement finie définie par
$U^\ast(\delta_y)=\sum_{x \in X} U(\delta_x)(y) \delta_x$ 
pour tout $y\in Y$.  
L'application $U^\ast$ est l'adjoint de $U$ pour les produit
naturels de $\CC(X,R)$ et $\CC(Y,R)$ et on a l'égalité $(U^\ast)^\ast=U$.

\subsection{Espaces propres généralisés}\label{sepgeneralises}
Soit $R$ un anneau, $M$ un $R$-module et $\Hecke$ un anneau agissant
sur $M$ par $R$-endomorphismes.

Pour tout caractère $\eta$ de $\Hecke$ à valeurs dans $R$, 
on définit le {\em sous-espace propre généralisé} $M(\eta)$ de $M$
pour $\eta$ comme le sous-module des vecteurs $m$ de $M$ tels que pour 
tout $T\in \Hecke$, 
il existe $n\geq 0$ vérifiant l'égalité $(T-\eta(T))^n(m) =0$. 
Si $R$ est un corps et $M$ est de dimension finie, 
il existe une extension finie $R'$ de $R$
et une décomposition $M\otimes R' = \oplus_\eta (M\otimes R')(\eta)$,
où $\eta$ décrit une famille finie de caractères de $\Hecke$ 
à valeurs dans $R'$.  

\section{Préliminaires locaux}\label{prelimslocaux}

Dans toute cette partie (sauf en~\ref{dec}), 
on fixe une place finie $v$ de $F$ 
inerte dans $E$. On note $w$ la place de $E$ au dessus de $v$ et 
$q$ le cardinal résiduel du corps $F_v$, complété de $F$ en $v$. 
Le groupe $G(F_v)$ est alors l'unique groupe unitaire à trois variables 
sur $F_v$ qui se déploie sur $E_w$.

Pour ne pas alourdir les notations, et dans cette partie uniquement, 
on note simplement $G$ le groupe $G_v=G(F_v)$, et $Z$ le centre de
$G$, qui est compact.
 
\subsection{Arbre} \label{arbre} Soit $\X$ l'immeuble de Bruhat-Tits  de $G$.
D'après~\cite{tits} ou \cite[1.4]{choucroun}, 
c'est un arbre, et  l'on a une décomposition 
de l'ensemble de ses sommets en deux parties $X \coprod X'$, 
tout sommet de $X$ (resp. $X'$) 
a $q^3+1$ (resp. $q+1$) voisins qui sont tous dans $X'$
(resp. $X$). Les points de $X$ sont les points {\it hyperspéciaux}
au sens de \cite{tits}, ceux de $X$ sont les points spéciaux qui ne
sont pas hyperspéciaux. On désigne par $A$ l'ensemble des arêtes
(non orientées) de $\X$. 

L'arbre $\X$ est muni d'une action par
automorphisme de $G$, le centre $Z$ agissant par l'identité. L'action
de $G$ sur $X$ (resp. $X'$) est transitive et le stabilisateur d'un sommet 
$x$ agit encore transitivement sur l'ensemble des 
sommets de $\X$ à distance $n$ de $x$ 
(\cite[1.4, 1.5]{choucroun}), et donc sur l'ensemble des 
éléments de $A$ d'origine $x$.

\subsection{Sous-groupes compacts}

\subsubsection{Compacts maximaux}
\label{compacts}

D'après~\cite{bt}, un sous-groupe compact maximal de $G$ fixe un
sommet de $\X$ et
un seul, ce qui définit une bijection entre l'ensemble des 
compacts maximaux de $G$ et $X \coprod X'$. 
Il y a donc deux classes de conjugaisons de
sous-groupes compacts maximaux de $G$, ceux qui fixent un sommet de
$X$, qu'on appelle {\it hyperspéciaux}, et ceux qui fixent un sommet
de $X'$, qu'on appelle {\it spéciaux}.

On suppose désormais donné $O \in X$ (resp. $O' \in X'$). 
On note $K$ (resp $K'$) le sous-groupe compact maximal de $G$ qui fixe 
 $O$ (resp $O'$), de sorte que l'on a 
les identifications canoniques $X = G/K=K \backslash G$ (resp. 
$X'= G/K'= K'  \backslash G$), la seconde étant induite
par l'application inverse de $G$.  

\subsubsection{Sous-groupe d'Iwahori}

On suppose aussi que $O$ et $O'$ sont voisins. Le stabilisateur 
$B=K \cap K'$ de l'arête $(O,O')$ est un sous-groupe d'Iwahori de $G$ 
et, par~\ref{arbre}, on a l'identification canonique $G/B=A$.

\subsection{Algèbres de Hecke}\label{operateurhecke}

\subsubsection{Algèbre de Hecke sphérique}
On note $\Hecke$
l'algèbre de convolution (pour la mesure de Haar $\mu_{K}$ de
volume $1$ sur $K$)
des fonctions sur $G$ à valeurs dans $\Z$, à
support compact, bi-$K$-invariantes
$$\Hecke =\CCC(K\backslash G /K,\Z).$$

Soit $T$ 
la fonction caractéristique dans 
$\CCC(G/K,\Z)=\CCC(X,\Z)$  
de l'ensemble des sommets à distance 2 de $O$. 
On a $T \in \Hecke$ 
et  les propriétés de transitivité de~\ref{arbre} impliquent 
l'égalité $\Hecke=\Z[T]$. 

\subsubsection{Algèbre de Hecke-Iwahori}
On appelle algèbre de Hecke-Iwahori le $\Z$-module 
muni du produit de convolution
$$\Hecke(G,B)=\CCC(B\backslash G /B,\Z)=\CCC(B\backslash A,\Z).$$
Soit $a\in \CCC(B\backslash A,\Z)$ (resp. $a'$) 
la caractéristique de l'ensemble des éléments de $A$ 
d'origine $O$ (resp. $O'$) et distincts de $(O,O')$. 
On définit $T_B \in \Hecke(G,B)$ par la formule
$$ T_B:= -a'a-aa'-(q^3-1)a'-(q-1)a-q^3(q-1).$$
On vérifie aisément que $a$ et $a'$ engendrent la $\C$-algèbre 
$\Hecke(G,B) \otimes \C$ et que $T_B$ est dans le centre de $\Hecke(G,B)$. 
\begin{remarque}
En fait, $T_B$ en est même un générateur (sur $\C$), mais nous n'utiliserons
pas ce fait. Ceci, ainsi que le lemme~\ref{relations}\,(vi), 
est d'ailleurs un petit fragment de la théorie du centre de 
Bernstein (\cite{bernstein}).
\end{remarque}

\subsection{Combinatoire des formes anciennes}\label{combinatoire}
  
\subsubsection{} 
Dans $\CCC(X,\Z)$, soit $U_1$ la fonction caractéristique
 de $O$, et $U_2$ la fonction caractéristique de l'ensemble des voisins
 de $O'$ distincts de $O$. Il est clair que $U_1$ et $U_2$ sont
 $B$-invariantes, i.e. appartiennent à 
$\CCC(B\backslash X,\Z)=\CCC(B \backslash G /K,\Z)$, 
qui est canoniquement muni d'une structure de $\Hecke$-module
par convolution à droite.  

\begin{prop}\label{lazarus}
Le $\Hecke\otimes_\Z \C$-module
$\CCC(B\backslash X,\C)$ 
est libre de base $(U_1,U_2)$.
\end{prop}
\begin{pf} D'après~\cite[prop. 7.2.4]{theselazarus},
le module $\CCC(B\backslash X,\C)$ 
est libre de rang $2$ sur $\Hecke\otimes \C$.
Pour conclure que $(U_1,U_2)$ est une base, il suffit de prouver que
cette famille est génératrice, ce qui se vérifie après
passage au quotient par tous les idéaux maximaux de 
$\Hecke \otimes \C=\C[T]$. Autrement dit, il suffit de montrer que pour tout
$\lambda\in \C$, $(U_1,U_2)$ est une base, c'est-à-dire une famille libre, 
du $\C$-espace vectoriel $[\CCC(X,\C)/(T-\lambda) \CCC(X,\C)]^B$ qui est 
de dimension $2$. 
Or, l'égalité $\lambda_1U_1+\lambda_2U_2 = (T-\lambda)f$, 
avec $\lambda_1,\lambda_2\in \C$ et $f \in \CCC(X,\C)$ 
implique $\lambda_1=\lambda_2=0$ par un argument de support.
\end{pf}
\begin{remarque} Dans~\cite[IV.4.8.2]{Joelthese}, 
cette proposition est démontrée plus généralement en rempla\c cant $\C$
par un anneau quelconque.
\end{remarque}
\begin{lemme}\label{deh} Si l'on voit $U_1$ et $U_2$ comme des fonctions de
 $\CCC(K \backslash G / B,\Z)=\CCC(K \backslash A,\Z)$, Alors $U_1$ 
est la fonction caractéristique de l'ensemble des arêtes passant 
par $O$,  et $U_2$ est l'ensemble des arêtes passant par un voisin de $O$, 
mais pas par $O$.
\end{lemme}
\begin{pf} 
C'est un calcul évident.
\end{pf}

\subsection{Opérateurs et relations}\label{operateurs}
Soit $R$ un anneau commutatif. 
Les fonctions a support fini $T \in \Hecke$, 
$T_B \in \Hecke(G,B)$, $U_1,U_2 \in 
\CCC(K\backslash G/B,\Z)$
induisent par convolution à droite des opérateurs localement finis 
$G$-équivariants à gauche entres les espaces
$\CC(X,R)$, $\CC(X',R)$ et $\CC(X',R)$
que l'on note par les mêmes lettres. Le lemme suivant est une traduction
immédiate des définitions.
\begin{lemme}
Notons $d$ la distance entre les sommets de $\X$. On a 
$$\begin{array}{crcl}
T \colon & \CC(X,R) &\rightarrow &  \CC(X,R) \\ 
 & \delta_x & \mapsto &  \ds{\sum_{d(y,x)=2}} \delta_y 
\end{array} $$
$$\begin{array}{crclccrcl}
U_1 : &\CC(X,R)& \rightarrow  &
\CC(A,R) & \  & U_2 : &\CC(X,R)& \rightarrow  &
\CC(A,R)  \\
& \delta_x & \mapsto & 
\ds{\sum_{d(x,x')=1}\delta_{(x,x')}}   & & & \delta_x &\mapsto &
\ds{\sum_{d(x',x)=1,\ d(y,x)=2}} \delta_{(x',y)}
\end{array}$$

$$ \text{et}\qquad  T_B= -a'a-aa'-(q^3-1)a'-(q-1)a-q^3(q-1) \Id,$$
où $a,a'\colon \CC(X,A) \to \CC(X,A)$ sont définis par les formules
$$ a\colon \delta_{(x,x')} \mapsto  \sum_{d(y',x)=1, y' \neq x'} 
\delta_{(x,y')}
\text{ \ \ et\ \ }
 a'\colon \delta_{(x,x')} \mapsto \sum_{d(y,x')=1, y \neq x} 
\delta_{(y,x')} $$
pour toute arête $(x,x')$ avec $x\in X$ et $x'\in X'$. 
\end{lemme}
\subsubsection{Définition}
On définit les opérateurs finis $G$-équivariants
$$ \begin{array}{crclccrcl} U : &\CC(X,R)&
\rightarrow &  
\CC(X',R) &\ \text{ et } \ & U' &: \CC(X',R) &\rightarrow   &\CC(X,R) \\ &\delta_x 
&\mapsto & \ds{\sum_{d(x',x)=1}} \delta_{x'}  & & &\delta_{x'}& 
\mapsto & \ds{\sum_{d(x',x)=1}} \delta_{x}.\end{array}$$

\begin{lemme}\label{relations} 
\begin{itemize} 
\item[(i)] $\ U'U = T + (q^3+1)$. 
\item[(ii)] L'opérateur $T$ 
est auto-adjoint et 
$U$ est l'adjoint de $U'$.
\item[(iii)] Soient $f_1,f_2 \in \CC(X,R)$. 
On a $U_1 f_1 + U_2 f_2=0$ 
si et seulement si il existe une constante $C\in R$ telle que $U f_2=C$ et 
$f_1-f_2=-C$.
\item[(iv)] Si $R$ est un corps, $U_1 f_1 + U_2 f_2=0$, et 
$T f_2 = \lambda f_2$, avec $\lambda \in R$ et $f_2 \neq 0$, 
on a $\lambda=q(q^3+1)$ ou bien $\lambda=-(q^3+1)$.
\item[(v)] Les opérateurs $U_1$ et $U_2$ sont injectifs.  
\item[(vi)] $U_i T =T_B U_i$ pour $i=1,2$. 
\end{itemize}
\end{lemme}
\begin{pf} Les preuves de 
(i), (ii) et (vi) sont des calculs évidents et laissées au lecteur
et (v) est triviale. 
Prouvons (iii).
Soit $(x,x')$ une arête de $A$, $x\in X$, $x'\in X'$. 
Le fait que 
$(U_1 f_1 + U_2 f_2)(\delta_{x',x})= 0$ s'écrit aussi 
$$f_1(x) + \sum_{y \neq x,\  y \text{ voisin de } x'} f_2(y) 
= 0,$$
soit $$f_1(x)-f_2(x) + (Uf_2)(x')=0.$$
Comme $\X$ est connexe, on en déduit qu'il existe 
$C\in R$ tel que $U f_2 = C$ sur $X'$,
$f_1-f_2 = -C$ sur $X$. 

Prouvons (iv). On a d'après (i) $U' U f_2 = U' (C) = (q^3+1) C =
 T f_2 + (q^3+1) f_2 = (\lambda + q^3+1) f_2$, d'où nous tirons
\begin{eqnarray} \label{eq} 
(q^3+1) C = (\lambda + q^3+1) f_2. 
\end{eqnarray}
Si $\lambda \not= -(q^3+1)$, alors la 
relation~(\ref{eq}) implique que $f_2$ est constante, d'où 
$Tf_2 = q(q^3+1) f_2$, et comme $f_2$ est non nul, 
il vient $\lambda = q(q^3+1)$. 
\end{pf}

\begin{prop}\label{ii} 
\begin{eqnarray*} 
\left( \begin{array}{cc}  U_1^\ast  U_1 & U_1^\ast  U_2 \\
U_2^\ast  U_1 & U_2^\ast  U_2 \end{array} \right) = 
 \left( \begin{array}{cc} q^3+1 & T \\ 
T & q(q^3+1)  + (q-1) T \end{array} \right)
\end{eqnarray*}
\end{prop}
\begin{pf} Pour toute arête $(x,x')$ de $A$, avec $x \in X$, 
l'adjoint  $U_1^\ast$ (resp. $U_2^\ast$) de $U_1$ (resp. $U_2$)
envoie $\delta_{(x,x')}$ 
sur $\delta_x$ (resp. sur $\sum_{y \neq x, 
\ d(y,x')=1} \delta_y$). La proposition est alors un calcul 
évident et laissé au lecteur.   
\end{pf} 

\begin{remarque} Le déterminant de cette matrice est un polynôme en $T$ dont 
les racines sont $-(q^3+1)$ et $q(q^3+1)$.
\end{remarque}

\subsection{Représentations non ramifiées}
\subsubsection{Notations}
Soit $R$ un anneau commutatif et $\pi$ une
représentation lisse à droite de $G$ sur $R$. 
Pour tout vecteur $u$ de l'espace de $\pi$ et toute fonction 
$f$ sur $G$ localement constante à
support compact, on définit 
$$ u \conv f:= \int_G f(g) (\pi(g)v) d\mu_K(g)$$
et pour tout sous-espace vectoriel $V$ de l'espace
de $\pi$, on pose $V \conv f :=\{u \conv f ;\ u \in V\}$.
Lorsqu'il n'y a pas d'ambiguïté, on note 
encore $f$ l'opérateur $u \mapsto u*f$. 
Ainsi, $T$ (
resp. $T_B$) 
désigne l'opérateur de convolution par $T$ (
resp. $T_B$) 
sur l'ensemble $\pi^{K}$ (
resp. $\pi^{B}$)  
des $K$-invariants (
resp. $B$-invariants) de $\pi$.
\begin{proposition} \label{oldu}
Soit $\pi$ une représentation complexe lisse de $G$ 
(éventuellement réductible) d'espace $V$, 
$\pi'$ la sous-représentation de $G$ engendrée 
par le sous-espace $\pi^{K}$ de $V$. Alors 
${\pi'}^{B} = \pi^{K} \conv U_1 + \pi^{K} \conv U_2$.
\end{proposition}
\begin{pf}
L'espace de ${\pi'}$ est l'image de l'application 
$$\pi^{K} \conv \CCC(K  \backslash G,\C)\to V,$$ 
et ${\pi'}^{B}$ est donc l'image de l'application 
$\pi^{K} \conv \CCC(K 
\backslash G/B,\C) \to V$. Le résultat découle donc de la 
proposition~\ref{lazarus}.
\end{pf}

\subsubsection{Tore maximal et sous-groupe de Borel}

On choisit un vecteur non nul $e_1 \in E_w^3$ isotrope pour la forme 
hermitienne standard, $e_2$ tel que $e_1$ et $e_2$ engendrent $(e_1)^\bot$, 
et enfin $e_3$ tel que $(e_1,e_2,e_3)$ soit une base. Dans
toute cette section, on décrit les éléments de $G$ 
par leur matrice dans la base $(e_1,e_2,e_3)$. Le tore diagonal 
$$D=\{\diag(a,b,c(a)^{-1})\ ; a \in E_w^\ast, b \in E_w^\ast, b c(b)=1\}$$
est un tore maximal de $G$, qui contient $Z=\{\diag(b,b,b)\ ;b c(b)=1\}$.
Soit $P$ le stabilisateur dans $G$ de la droite $(e_1)$ : c'est un 
sous-groupe de Borel contenant $D$.
 
\subsubsection{Induites non ramifiées}

Pour tout 
$\alpha \in \C$, on définit $I(\alpha)$ comme la représentation 
induite unitaire de $P$ à $G$ du caractère de $D$~:
$$\diag(a,b,c(a)^{-1}) \mapsto \alpha^{\val(a)}.$$
Comme ce caractère est trivial sur $Z$, $I(\alpha)$ est de caractère central 
trivial.

\subsubsection{Structure des induites non ramifiées}
\label{induites}
Soit $\alpha\in \C$; 
d'après~\cite[page 126]{KE} ou \cite[th. 2.4.6 et prop 2.4.7]{choucroun},
la représentation $I(\alpha)$ est indécomposable et vérifie 
les propriétés suivantes. 
\begin{itemize} 
\item[(1)] Si $\alpha \not=q^{\pm2}$ et $\alpha \not= -q^{\pm1}$, 
alors $I(\alpha)$ est irréductible.
\item[(2)] Si $\alpha = q^{\pm 2}$, alors $I(\alpha)$ a deux 
facteurs de Jordan-Hölder, la Steinberg $\St$ et la représentation unité.
\item[(3)] Si $\alpha = -q^{\pm 1}$,
 alors $I(\alpha)$ a deux facteurs de Jordan-Hölder, 
l'un non ramifié noté
$\pi^n$ et l'autre ramifié noté $\pi^s$.
D'après~\cite[page 396]{rog3}, la représentation $\pi^n$ 
est non tempérée et la représentation $\pi^s$ est de 
carré intégrable. 
\end{itemize}
\subsubsection{Invariants des induites non ramifiées}
\label{invariants}
\begin{itemize}
\item[(1)] Une représentation admissible irréductible possède 
un vecteur  $B$-invariant
si et seulement si c'est un facteur de $I(\alpha)$ pour 
un certain $\alpha\in \C$
(\cite{cartier}).
\item[(2)] Pour tout $\alpha \in \C$, $\dim I(\alpha)^{K}=\dim
I(\alpha)^{K'}=1$ et $\dim I(\alpha)^B=2$ (\cite{cartier} ou
\cite[IV.4.8.2 et IV.6.2.5]{Joelthese}).
\item[(3)] $\St^{K}=\dim \St^{K'} = 0$, $\dim (\pi^n)^{K}=
\dim (\pi^s)^{K'}=1$, $(\pi^n)^{K'} = (\pi^s)^{K}=0$ 
(\cite[prop. 2.4.7]{choucroun} ou \cite[IV.4.8.2 et IV.6.2.5]{Joelthese}). 
\end{itemize}

\subsubsection{Matrices de Hecke}
\label{heckeinerte}
Pour tout $\alpha\in \C$, 
définissons la matrice de Hecke $\mat_{\pi}$ (à conjugaison près)
de l'unique sous-quotient non ramifié $\pi$ de 
$I(\alpha)$ comme $$\mat_\pi = \diag(\alpha,1,\alpha^{-1}) \in \Sl_3(\C).$$
Ainsi $\mat_{\pi^n} = \diag(-q,1,-q^{-1})$.

Par exception, dans la théorie globale (partie~\ref{sectionrepresentations}), 
on note $\mat_{\pi,w}$ pour 
$\mat_{\pi}$ au lieu de $\mat_{\pi,v}$ : cela tient à ce que, à strictement 
parler, $\mat_\pi$ est la matrice de Hecke du changement de base de $\pi$ à 
$E_w$ et non celle de $\pi$.

\subsection{Transformation de Satake}
\label{satakeinerte}
\begin{lemme}
Soit $\alpha\in \C$; l'opérateur $T$ agit sur la droite 
$I(\alpha)^{K}$ par l'homothétie de rapport
$$q^2 (\alpha + \alpha^{-1}) + q-1.$$
En particulier, si $\alpha=q^{\pm 2}$ (resp. $\alpha=-q^{\pm 1}$), alors 
la valeur propre de $T$ est $q(q^3+1)$ (resp. $-(q^3+1)$). 
\end{lemme}
\begin{pf} Voir~\cite[prop 3.1.2]{choucroun}. \end{pf}
\begin{lemme}\label{vptb}
L'opérateur $T_B$ agit sur $\st^B$ (resp. $(\pi^s)^B$) 
par l'homothétie de rapport $q(q^3+1)$ (resp. $-(q^3+1)$).
\end{lemme}
\begin{pf}
Comme $I(\alpha)$ est indécomposable (cf.~\ref{induites}), il 
découle de~\cite{borel2} que $I(\alpha)^B$
est un module indécomposable sur $H(G,B)\otimes \C$, donc $T_B$ y agit
par une homothétie. En particulier, par le lemme~\ref{relations} (vi), 
le rapport de $T_B$ agissant sur $I(\alpha)^B$ est égal à celui de 
$T$ agissant sur $I(\alpha)^K$.
Le lemme~\ref{vptb} découle donc du lemme précédent. 
\end{pf}

\subsection{Théorie locale en une place décomposée}
\label{dec}
\label{heckedec}
Soit $v$ une place décomposée de $F$. Le choix d'une place $w$ de $E$
divisant $v$ détermine (à isomorphisme intérieur près) 
un isomorphisme $G_v \simeq \Gl_3(E_w)=\Gl_3(F_v)$. Si $\pi$ est une
représentation non ramifiée de $G_v$, on note $\mat_{\pi,w}$ la matrice
de Hecke de $\pi$ vue comme représentation de $\Gl_3(E_w)$ via
l'isomorphisme précédent. Si $\bar w$ est l'autre place de $E$
au-dessus de $v$, on a, à conjugaison près, $\mat_{\pi,w} =
\mat_{\pi,\bar w}^{-1}$.
  
\section{Sorites sur les congruences} 

Dans cette section on rappelle pour le confort du lecteur quelques résultats 
faciles, bien connus et souvent utilisés dans le contexte des congruences 
entre formes automorphes (voir les travaux de Hida, Ribet, etc.), et
on les étend quelque peu.

Dans cette partie, $\anneau$ désigne un 
produit fini d'anneaux de valuation discrète
et $L$ (resp. $\muz$) 
le produit des corps de fractions (resp. des idéaux maximaux) 
de ces anneaux. On désigne par $M$ un $\anneau$-module libre de type fini.

\subsection{Sous-modules saturés}


Pour tout sous-module $N$ de $M$, on pose 
$N^\sat:=(N \otimes L) \cap M$ dans $M \otimes L$. 
On a $N=N^{\sat}$ si et 
seulement si $N$ est facteur direct de $M$.

\begin{lemme}\label{noyau}
Soit $u\colon  N \rightarrow M$ un morphisme injectif de 
$\anneau$-module libre de rang fini.
Pour tout entier $\alpha$ assez grand, 
$u(N)^\sat/u(N)$ et $\Ker(u\otimes(\anneau/\muz^\alpha))$ sont
isomorphes en tant 
que $\anneau$-modules. 
\end{lemme}
\begin{pf}
On peut supposer que $\anneau$ est un anneau de
valuation discrète et que $u(N)^\sat=M$. Par le théorème 
des diviseurs élémentaires, nous sommes ramenés 
à ne traiter que le cas évident où $M=\anneau$. 
\end{pf}

\subsection{Congruences}

Jusqu'à la fin de cette partie, $A$ et $B$ désignent des sous modules de $M$
vérifiant $A=A^\sat$, $B=B^\sat$ et $A \cap B=0$, 
$\mu$ désigne un idéal maximal de $\anneau$ et $c \geq 1$ un entier.  

\subsubsection{}
On dit qu'il existe une 
{\it congruence modulo $\mu^c$} entre $A$ et $B$ dans $M$, s'il existe 
$f \in A \setminus\mu A$, $g \in B\setminus \mu B$, 
avec $f-g \in \mu^c M$. On appelle
{\it module de congruence} entre $A$ et $B$ (dans $M$) 
le $\anneau$-module $(A \oplus B)^\sat / (A \oplus B)$. 
Ce module tient son nom du lemme suivant.

\begin{lemme} \label{modcong} Il existe des congruences  
modulo $\mu^c$ entre $A$ et $B$
dans $M$ si et seulement si le module de 
congruence contient un sous-module isomorphe à $\anneau/\mu^c$.
\end{lemme}
\begin{pf}
Soit $\pi$ un générateur de $\mu$.
Si le module de congruence contient un sous-module isomorphe à $\anneau/\mu^c$ 
alors il existe $h \in (A \oplus B)^\sat$, avec 
$\pi^c h \in A \oplus B$ mais $\pi^{c-1} h \not \in A \oplus B$. 
On a alors $\pi^c h = f +g$ avec $f \in A$ et $g \in B$, et
puisque $B=B^\sat$, on a $f \in \pi A$ si et seulement si $g  \in \pi B$, 
ce qui n'est pas car $\pi^{c-1}h \not\in A + B$. 
On a donc une congruence modulo $\pi^c$ entre $f$ et $g$. 

Réciproquement, si $f$ et $g$ sont en congruences modulo $\pi^c$, alors 
l'élément
$h:=(f-g)/\pi^c \in (A \oplus B)^{\sat}$ est tel que 
$\pi^{c-1} h \not \in (A \oplus B)$. 
\end{pf}

\subsection{Produit hermitien et congruences}

\label{phc} 
Soit $\gamma\colon \anneau \to \anneau$ 
un automorphisme involutif d'anneaux. 
Pour tout module $N$ libre de type fini sur $\anneau$, 
on note $N^\anti$ le $\anneau$ module des formes $\gamma$-linéaires
(i.e. telles que $f(a n)=\gamma(a) f(n)$ pour $a \in \anneau$, $n \in N$)
sur $N$ et on appelle {\it produit hermitien} sur $N$ tout 
morphisme $\anneau$-linéaire 
$p_N: N \rightarrow N^\anti$ tel que $p_N(x)(y)=\gamma(p_N(y)(x))$ 
pour tout $x,y \in N$; on dit que $p_N$ est {\it non dégénéré}
 s'il est bijectif. 
Pour toute partie $P$ de $N$, on note alors $P^\bot$ l'ensemble des $x\in N$ 
tel que $p_N(P)(x)=0$ ; c'est un sous-module facteur direct de $N$. 
Si $P$ est un sous-module de $N$, on note 
$p_P : P \rightarrow P^\anti$ la restriction de $p_N$ à $P$. 

Jusqu'à la fin de cette partie, on désigne par $p_M$ un  
produit hermitien non dégénéré sur $M$ tel que $B \subset A^\perp$. 

\begin{lemme} \label{modprod}
On suppose que $(A+B)^\bot \oplus A$ est facteur direct dans $M$.
Alors le module de congruence entre $A$ et $B$ dans $M$ est isomorphe à 
$A^\anti/p_A(A)$.
\end{lemme}
\begin{pf}
Le morphisme naturel de $\anneau$-modules 
$(A\oplus B)^\sat/(A\oplus B) \to A^\anti/p_A(A)$
induit par $p_M$
est injectif car $B$ est orthogonal à $A$ et surjectif
car toute forme antilinéaire sur $A$ se prolonge par
hypothèse en une forme antilinéaire sur $M$ nulle sur $(A\oplus B)^\bot$, 
qui du fait de la non-dégénérescence de $p_M$ est l'image par 
$p_M$ d'un élément de $(A\oplus B)^\sat$. 
\end{pf}

Si $u:\,N \rightarrow M$ est un morphisme injectif 
de modules libres munis de produit
scalaires $p_N$ et $p_M$ non dégénérés, on note 
${}^t u \colon M^\anti \to N^\anti$ la transposée de $u$
et on appelle {\it adjoint} de $u$ le 
morphisme $u^\ast=p_{N}^{-1} \circ {}^t u \circ p_M$.
\begin{lemme} \label{longueur}
Pour tout entier $\alpha$ et tout idéal maximal $\mu$ de $\anneau$, on a
$$\lng\Big(\ker (u \otimes(\anneau/\mu^\alpha))\Big) 
\leq \frac12{\val_\mu(\det(u^\ast u))}.$$
\end{lemme}
\begin{pf} 
Soit $P=u(N)^\sat$. Par le lemme~\ref{noyau}, la 
longueur à calculer est majorée par la $\mu$-valuation du 
déterminant de $u_P\colon N \to P$ dans deux bases
quelconques de $N$ et $P$. Or, puisque $p_N$ est non dégénéré, 
$u^*_P \colon P \to N$ est bien défini et $u^*u=u_P^* u_P$. 
Puisque $\det(u_P)$ divise $\det(u^*_P)$ dans $\anneau$, le lemme est prouvé. 
\end{pf} 

\section{Augmentation du niveau}

\subsection{Espaces de formes automorphes}\label{espacesdeformes}

\subsubsection{Niveau, type et poids}\label{ntp}

On fixe un sous-groupe compact ouvert $K=\prod_v K_v$ de $G(\A_{F,f})$
({\it le niveau}). On note $\Sigma$ l'ensemble fini (\cite{tits}) 
des places $v$ où $K_v$ n'est pas un compact
maximal hyperspécial et pour toute place $\mu$ d'un corps de nombres $L$,  
on désigne par $\Sigma_\mu$ la réunion de $\Sigma$ et
des places finies $v$ de $F$ de même caractéristique résiduelle que $\mu$.
Pour tout ensemble $\Sigma'$ de places de $F$ contenant $\Sigma$ 
on définit l'algèbre de Hecke 
$$\Hecke^{\Sigma'} :=  \Hecke(\textstyle{\prod_{v \not \in S}' G_v,
\prod_{v \not \in S} K_v})$$
qui est aussi le produit tensoriel restreint des $\Hecke(G_v,K_v)$ pour
$v\not\in \Sigma'$.
Cette algèbre contient l'opérateur de Hecke $T_v$
défini en~\ref{operateurhecke}, 
pour toute place $v \not\in \Sigma'$ inerte dans $E$.

On fixe aussi une représentation complexe 
lisse irréductible $(J,V)$ de $\prod_{v \in \Sigma} K_v$ ({\it le type}),
que l'on voit indifférement comme une représentation de $K$.
Enfin, on fixe une représentation complexe continue irréductible
$(\rho,W)$ de $G(\A_{F,\infty})$ ({\it le poids}).  
On note $\rho^\ast$ (resp. $J^\ast$) la représentation duale de $\rho$
(resp. $J$).

\subsubsection{Espaces de formes automorphes complexes}
\label{parformescomplexes}\label{defB}
On appelle espace des {\em formes automorphes complexes}
de niveau $K$, de type $J$ et de poids $\rho$  
l'espace vectoriel
$$S_{K,J,\rho,\C}:= \Hom_{K \times G(\A_{F,\infty})}(J \otimes \rho,\BB),$$
où $\BB$ désigne 
l'espace des fonctions complexes lisses (i.e. invariante par
translation par un ouvert de $G(\A_{F,f})$ et $\CC^\infty$ sur 
$G(\A_{F,\infty})$) sur $G(F)\backslash
G(\A_{F})$. 
De manière trivialement équivalente, on peut voir 
$S_{K,J,\rho,\C}$ comme l'espace des fonctions
complexes lisses $f$ sur $G(F)
\backslash G(\A_F)$ à  valeur dans  
$V^\ast \otimes W^\ast$ et vérifiant 
\begin{eqnarray}
\label{eqfonc1} && 
f(gku_{\infty})=(J^\ast(k)^{-1}\otimes \rho^\ast(u_\infty)^{-1})
f(g) \text{ pour $g \in G(\A_{F}),\ k \in K,\ 
u_{\infty} \in G(\A_{F,\infty})$}.  
\end{eqnarray}
On peut ainsi voir $S_{K,J,\rho,\C}$ comme le sous-espace 
des $\prod_{v \neq \Sigma} K_v$-invariants
d'un espace de fonctions complexes sur $G(F)
\backslash G(\A_F)$ sur lequel 
$G(\A_F)$ agit à gauche par translation à droite,
de sorte que $S_{K,J,\rho,\C}$ est muni d'une structure de
$\Hecke^\Sigma$-module à gauche.

Si $f$ vérifie
l'équation fonctionnelle~(\ref{eqfonc1}), la restriction $f'$ de
$f$ à $G(\A_{F,f})$ vérifie 
\begin{eqnarray} \label{eqfonc2} && 
f'(ugk)=(J^\ast(k)^{-1}\otimes \rho^\ast(u_\infty))
f(g) \text{ pour tout $g \in G(\A_{F,f}),\ u \in G(F),\ k \in K$}.  
\end{eqnarray}
et détermine $f$ par densité.

\subsubsection{Modèle sur un corps de nombres}\label{modelecorpsdenombres}

La représentation $(J,V)$ de $K$ se factorise par un quotient fini. 
Il existe donc un corps de nombres $L \subset \C$
et un $L$-modèle $(J_L,V_L)$ de $(J,V)$.
Quitte à grossir $L$, on peut supposer qu'il est stable par la
conjugaison complexe, ce qui donne un sens à la notion 
de produit hermitien sur $V_L$. Utilisant que $J_L(K)$ est fini, 
nous choisissons un produit hermitien sur $V_L$ stabyle par $J_L$.

Soit $\rho_\C$ le prolongement de 
la représentation $\rho$  du groupe alg\'ebrique r\'eel
$G(\A_{F,\infty})=\prod_{\sigma|\infty} G(\R)$ sur $W$ en
une représentation de $\prod_{\sigma|\infty} G(\C)$ 
sur le même espace $W$. Comme $\rho_\C$ est algébrique, elle   
admet un modèle sur un corps de nombres. 
Quitte à grossir  $L$, on peut supposer qu'il existe un 
tel $L$-modèle $(\rho_L,W_L)$, $\rho_L$ étant une représentation
de $\prod_{\sigma|\infty} G(\sigma\colon F \to L)$, 
que l'extension $L/\Q$ 
est galoisienne et qu'enfin $L$ contient $F$ et n'est pas inclus dans $\R$. 
 
Notons $L_0$ le sous-corps $L \cap \R$ de $L$, de sorte que $[L:L_0]=2$.
Comme $G(\A_{F,\infty})$ est compact, 
$\rho$ laisse stable un produit hermitien sur $W$, 
et il en va donc de même de la restriction de $\rho_L$ à 
$\prod_{\sigma| \infty}G(L_0)$. Nous munissons 
$W_L$ d'une telle structure hermitienne. 

Pour tout $u\in G(F)$, et tout plongement $\sigma$ de $F$ dans $\C$,
on a $\sigma(u) \in L_0$ ; il s'ensuit que 
l'élément $\sigma^\ast(u_\infty)$ de $\Gl(W)$ est en fait un élément unitaire
de $\Gl(W_L)$.
Par ailleurs, \cite[5.1]{borel} assure que
$G(F)\backslash G(\A_{F,f})$ est compact, et donc est
muni d'une unique mesure de probabilité invariante à droite.
Ceci donne un sens à la définition suivante. 
\begin{definition}
Soit $S_{K,J,\rho,L}$ 
le $L$-espace vectoriel des fonctions lisses de $G(\A_{F,f})$ dans 
$V_L^\ast \otimes W_L^\ast$ vérifiant l'équation
fonctionnelle~(\ref{eqfonc2}), muni du produit hermitien 
obtenu en intégrant celui de $V_L^\ast \otimes W_L^\ast$
sur $G(F)\backslash G(\A_{F,f})$.
Pour toute $L$-algèbre $R$, on définit le $R$-module $S_{K,J,\rho,R}$
des formes automorphes à valeurs dans $L$ par 
$$S_{K,J,\rho,R} := S_{K,J,\rho,L}\otimes_L R,$$ 
c'est un $R$-module hermitien si $R$ est muni d'une involution 
prolongeant la conjugaison complexe de $L$.
\end{definition}
Procédant comme en~\ref{parformescomplexes}, on peut
voir $S_{K,J,\rho,R}$ comme le sous-espace des $\prod_{v \neq \Sigma} K_v$
invariants d'une représentation unitaire de $\prod'_{v \neq \Sigma} G_v$. 
Puisque $\Hecke^\Sigma$ est commutative, ses éléments 
agissent normalement sur $S_{K,J,\rho,R}$, qui
est donc un $\Hecke^\Sigma \otimes_\Z R$-module 
semi-simple lorsque $R$ est un corps. 

\subsubsection{}\label{pardescription} Explicitons maintenant cette définition.
Puisque $G(F)\backslash G(\A_{F,f})$ est compact,
il existe des éléments $x_1,\ldots,x_h$ de $G(\A_{F,f})$ 
tels que 
$$G(\A_{F,f})=\coprod_{i=1}^h G(F)x_i K.$$
Pour tout $i=1,\ldots,h$, 
le groupe $\Gamma_i=(x_i^{-1} G(F) x_i) \cap K$ 
est fini car compact et discret, et on a un isomorphisme 
\begin{eqnarray} \label{description} 
S_{K,J,\rho,L} &\rightarrow &
\oplus_{i=1}^h (V_L^\ast \otimes W_L^\ast)^{\Gamma_i}\\
f &\mapsto & (f(x_1),\dots,f(x_h)).\nonumber \end{eqnarray}
D'autre part, le produit hermitien
sur $S_{K,J,\rho,L}$ est à une constante rationnelle près 
la somme sur $i=1,\ldots,h$ de ceux des espaces 
$(V_L^\ast \otimes W_L^\ast)^{\Gamma_i}$, 
pondérée par les inverses des cardinaux des sous-groupes $\Gamma_i$.

On a montré l'énoncé suivant.
\begin{lemme}
L'espace $S_{K,J,\rho,L}$ 
est un $\Hecke^\Sigma$-modèle de $S_{K,J,\rho,\C}$. 
\end{lemme}

\subsubsection{Décomposition de $S_{K,J,\rho,L}$}
\label{decomp} 
Appliquant~\ref{sepgeneralises} au $\Hecke^\Sigma\otimes L$-module
$S_{K,J,\rho,L}$, on peut supposer, 
quitte à remplacer $L$ par une extension 
finie, qu'il existe un ensemble fini $E_{K,J,\rho,L}$ de caractères de
$\Hecke^\Sigma$ à valeurs dans $L$ et une décomposition
\begin{eqnarray*}
S_{K,J,\rho,L}=\oplus_\eta S_{K,J,\rho,L}(\eta)
\end{eqnarray*}
(en fait, puisque $S_{K,J,\rho,L}$ est un $\Hecke^\Sigma\otimes L$-module 
semi-simple, les sous-espaces propres généralisés $S_{K,J,\rho,L}(\eta)$
sont de simples sous-espaces propres pour l'action de $\Hecke^\Sigma$). 

En particulier, dans le cas où $J=1_\C$ et $\rho=1_\C$ sont
les actions triviales sur $\C$, 
la droite des fonctions constantes
est un sous-$\Hecke^\Sigma\otimes_\Z L$-module 
de $S_{K,1_\C,1_\C,L}$
sur lequel $\Hecke^\Sigma$ agit par un caractère à valeurs dans $\Z$ 
noté $\eta_\const$. Pour toute place inerte $v \not \in \Sigma$,
on a $\eta_\const(T_v)=q(q^3+1)$, où $q$ désigne le cardinal résiduel de $F_v$.

Rappelons qu'on dit que deux éléments $x$ et $y$ de $L$
sont congrus modulo $\mu$ si $x-y \in \mu \anneau_L$, où 
$\anneau_L$ désigne l'anneau des entiers de $L$. 
Le lemme suivant est utile pour séparer les caractères de $\Hecke^\Sigma$. 
\begin{lemme}\label{egalitecarac}
Soient $\eta$, $\eta'$ deux caractères de $E_{K,J,\rho,L}$. 
Si pour presque toute place $\mu$ de $L$, les restrictions de
$\eta$ et $\eta'$ à $\Hecke(G_v,K_v)$ sont congrues modulo $\mu$,
alors $\eta$ et $\eta'$ sont congrus modulo $\mu$. 
\end{lemme}
\begin{pf}
Cela résulte immédiatement de l'existence des représentations
galoisiennes de $\Gal(\bar E/E)$ attachées à $\eta$ et à $\eta'$ 
(cf.~\ref{repattachee}) et du théorème de $\check{\text{C}}$ebotarev.
\end{pf}
\begin{lemme}\label{defs}
 Il existe un ensemble fini $S_1$ de nombres premiers (ne dépendant que
de $(K,J,\rho)$) tel que pour toute
place $\mu$ de $L$ ne divisant aucun élément de $S_1$
et tout couple de caractères $\eta$, $\eta'$ 
de $E_{K,J,\rho,L} \cup \{\eta_\const\}$, 
$\eta$ est congru à $\eta'$ si et seulement si $\eta=\eta'$.
\end{lemme}
\begin{pf} C'est évident. 
\end{pf}

\subsection{Modèles entiers}\label{modeleentier} 
Nous n'allons pas définir de structures entières ``globales'' des
espaces précédents, mais seulement en presque toute place $\mu$ 
de $L$. 
\subsubsection{Définition de $\rho_\muz$.}
Soit $\muz$ une place finie de $L_0=L\cap\R$. 
Notons $L_\muz$ le complété de $L$ en $\muz$ et
$\anneau_\muz$ son anneau d'entiers, qui est 
une algèbre semi-locale sur l'anneau des entiers du corps local 
$(L_0)_\muz$, munie d'une unique involution prolongeant
celle de $\anneau_L$. 

Pour toute place infinie $\sigma$ de $F$ correspondant à un plongement 
$F \to L_0$, notons $\sigma^{-1}(\muz)$
la trace de $\muz$ sur $F$ par ce plongement, de
sorte que $F_{\sigma^{-1}(\muz)}$ est plongé dans $L_\muz$.  
On a donc un morphisme naturel de groupes, 
compatible aux plongements de $G(F)$,
\begin{eqnarray} \label{morphismemu}
 &&  G(\A_{F,f}) \to \prod_{\sigma|\infty} G(F_{\sigma^{-1}(\muz)}) \to 
\prod_{\sigma|\infty} G(L_\muz). 
\end{eqnarray}
Définissons une représentation $\rho_\muz$ de $G(\A_{F,f})$ 
sur $W^*_{L_\muz}$ 
en composant (\ref{morphismemu}) avec 
$$\rho_{L_\muz} \colon \prod_{\sigma|\infty} 
G(L_\muz)\to \GL(W^\ast_{L_\muz}).$$
Pour tout $u\in F$, l'élément $\rho_\muz(u)$
de $\GL(W^\ast_{L_\muz})$ est en fait dans $\GL(W^\ast_L)$ 
et on a l'égalité
\begin{eqnarray}\label{egaliteu}  
\rho_\muz(u)&=&\rho(u_\infty).  
\end{eqnarray}

Enfin, on fixe un $\anneau_L$-réseau $V^\ast_{\anneau_L}$ 
(resp. $W^\ast_{\anneau_L}$) de $V^\ast_{L}$ (resp. $W^\ast_L$).

\subsubsection{Définition des modèles entiers}

\begin{definition}
Pour toute place $\muz$ de $L_0$, soit 
$S_{K,J,\rho,\anneau_\muz}$
le sous $\anneau_\muz$-module de 
$S_{K,J,\rho,L_\muz}$ 
constitué des fonctions $f$ telles que
$$  \rho_\muz^\ast(g)^{-1}f(g) 
\in V_{\anneau_\muz}^\ast \otimes W_{\anneau_\muz}^\ast 
\hbox{ pour tout $g \in G(\A_{F,f})$.}$$
\end{definition}

Par définition, $\rho_\muz$ se factorise par 
$\prod_{v} G(F_v)$, où $v$ décrit l'ensemble des places de $F$
de même caractéristique résiduelle que $\muz$, 
de sorte que l'on 
peut donc encore voir $S_{K,J,\rho,\anneau_\muz}$ comme l'ensemble
des $\prod_{v \not\in \Sigma_\muz} K_v$-invariants 
d'un espace de fonctions sur $G(\A_{F,f})$ sur lequel 
$\prod'_{v \not\in \Sigma_\muz} G_v$ agit par translation à droite.  
On en déduit que l'action de $\Hecke^{\Sigma_\muz}$ sur  
$S_{K,J,\rho,L_\muz}$ laisse stable 
$S_{K,J,\rho,\anneau_\muz}$.

\subsubsection{D\'efinition de $S_2$.}\label{modelegroupe}
Le morphisme de $L$-schémas en groupes
$$\prod_{\sigma|\infty} G \times_\sigma \spec L  
\to \mathop{\mathrm{\bf GL}}(V_L^\ast \otimes W_L^\ast)$$ 
induisant $\rho_L$ se prolonge en un morphisme 
$$\prod_{\sigma|\infty} G \times_\sigma \spec \anneau_L 
\to \mathop{\mathrm{\bf GL}}(V_{\anneau_L}^\ast \otimes W_{\anneau_L}^\ast)$$ 
au dessus d'un ouvert de $\spec \anneau_L$, d'où l'existence 
d'un ensemble fini $S_2$ de nombres premiers, ne
dépendant que de $(K,J,\rho)$, tel que pour toute place
finie $\muz$ de $L\cap\R$ ne divisant aucun élément de $S_2$, 
le réseau $V_{\anneau_\muz}^\ast \otimes W^\ast_{\anneau_\muz}$ 
de $V^\ast_{L_\muz}\otimes W^\ast_{L_\muz}$ 
est stable sous la restriction de
$J_{L_\muz}^\ast \otimes \rho_{L_\muz}^\ast$ 
au sous-groupe $\prod_{\sigma|\infty} G(\anneau_\muz)$
de $\prod_{\sigma|\infty} G(L_\muz)$.

\subsubsection{D\'efinition de $S_3$.}\label{stabilitereseau}
Par définition de la topologie adélique, il existe  un ensemble $S_3$
de nombres premiers tel que pour toute place finie $v$ de $F$ ne
divisant aucun élément de $S_3$, on a l'égalité $K_v=G(\anneau_v)$.
Pour toute place finie $\muz$ de $L_0$ ne divisant
aucun élément de $S_2\cup S_3$, 
\ref{modelegroupe} montre alors que la restriction
de $\rho_\muz \otimes J_L$ à $K$ est à valeurs dans  
$\GL(V_{\anneau_L}^\ast\otimes W_{\anneau_L}^\ast)$. 

\subsubsection{D\'efinition de $S_4$.}\label{produitscalaire}
Comme en~\ref{pardescription}, considérons une famille 
$x_1,\ldots,x_h$ de représen\-tants de $G(F)\backslash G(\A_{F,f})/K$ 
et, pour $i=1,\ldots,h$, posons $$\Gamma_i= (x_i^{-1}G(F) x_i)\cap K.$$ 
Le sous-$\anneau_L$-module 
$(V_{\anneau_L}^\ast \otimes W_{\anneau_L}^\ast)^{\Gamma_i}$
du $L$-espace vectoriel 
$(V_{L}^\ast \otimes W_{L}^\ast)^{\Gamma_i}$ en est un réseau. 
De plus, il existe un ensemble fini de nombres premiers $S_4$, 
tel que pour toute place $\muz$ de $L\cap \R$ 
ne divisant aucun élément de $S_4$
et pour tout $i=1,\ldots,h$, la restriction du produit scalaire
de $V^\ast_{L_\muz}\otimes W^\ast_{L_\muz}$
à
$(V^\ast_{\anneau_\muz}\otimes W^\ast_{\anneau_\muz})^{\Gamma_i}$
est à valeurs dans $\anneau_\muz$ et non-dégénérée. 
Quitte à remplacer $S_4$ par un ensemble le contenant, on 
suppose de plus que $S_4$ ne contient aucun des
diviseurs des cardinaux des groupes finis $\Gamma_1,\ldots,\Gamma_h$. 
\begin{lemme}\label{lemmemodeleentier}
Soit $\muz$ une place finie de $L_0$ 
ne divisant aucun élément de $S_2\cup S_3\cup S_4$, 
alors $S_{K,J,\rho,\anneau_\muz}$ est un 
$\Hecke^{\Sigma_\muz}$-modèle 
de $S_{K,J,\rho,L_\muz}$, sur lequel le produit hermitien est non dégénéré.
\end{lemme}
\begin{pf} 
Montrons que $S_{K,J,\rho,\anneau_\muz}$ est un réseau de 
$S_{K,J,\rho,L_\muz}$.
Pour tout $f \in S_{K,J,\rho,\anneau_\muz}$, 
$g \in G(\A_{F,f})$, $u\in G(F)$ et $k\in K$, on a par~(\ref{eqfonc2}) 
et~(\ref{egaliteu})
\begin{eqnarray*}
  \rho_\muz^\ast(ugk)^{-1} f(ugk)
&= &(J^\ast(k)^{-1} \otimes\rho_\muz^\ast(k)^{-1}\rho_\muz^\ast(g)^{-1}
\rho_\muz^\ast(u)^{-1}\rho^\ast(u_\infty))f(g)\\
&= &(J^\ast(k)^{-1} \otimes\rho_\muz^\ast(k)^{-1})
(\rho_\muz^\ast(g)^{-1}f(g)).\end{eqnarray*}
Sous les notations de~\ref{produitscalaire}, 
le calcul précédent et \ref{stabilitereseau} impliquent que 
$f\in S_{K,J,\rho,\anneau_\muz}$ si et seulement si
$\rho_\muz^\ast(x_i)^{-1}f(x_i)\in 
V_{\anneau_\muz} \otimes W_{\anneau_\muz}$ pour $i=1,\ldots,h$, et donc 
l'isomorphisme~(\ref{description}), avec $R=L_\muz$, 
induit un isomorphisme  
$$ S_{K,J,\rho,\anneau_\muz} \simeq 
\bigoplus_{i=1}^h \Big(\rho_\muz^\ast(x_i)
(V_{\anneau_\muz} ^\ast \otimes 
W_{\anneau_\muz}^\ast)\Big)^{\Gamma_i}.$$
Or, pour $i=1,\ldots,h$, 
le $\anneau_\muz$-module $\Big(\rho_\muz^\ast(x_i)
(V_{\anneau_\muz} ^\ast \otimes  
 W_{\anneau_\muz}^\ast)\Big)^{\Gamma_i}$  
 est un réseau 
 du $L_\muz$-espace vectoriel 
$(V_{L_\muz} ^\ast \otimes W_{L_\muz}^\ast)^{\Gamma_i}$.  

Il reste à prouver que 
la restriction du produit hermitien $S_{K,J,\rho,\anneau_\muz}$
à $S_{K,J,\rho,\anneau_\muz}$ est à valeurs dans $\anneau_\muz$
et non dégénérée, ce qui découle de~\ref{produitscalaire} 
puisque ce produit hermitien est  
une somme directe, pondérée par des éléments inversibles de $\anneau_\muz$,
de produits hermitiens à valeurs dans $\anneau_\muz$ et non dégénérés.
\end{pf}

\subsubsection{Définition} On pose 
$S_{K,J,\rho,R} : = S_{K,J,\rho,\anneau_\muz} \otimes_{\anneau_\muz} R$
pour toute $\anneau_\muz$-algèbre $R$. Cette définition est compatible
avec~\ref{modelecorpsdenombres} par le lemme précédent.

 \begin{remarque}
Le lemme~\ref{lemmemodeleentier} montre en particulier que
la restriction à $\Hecke^{\Sigma_\muz}$ des 
caractères de $\Hecke^\Sigma$ définis en~\ref{decomp}
est à valeurs dans l'anneau localisé en $\mu$ de $\anneau_L$
pour presque toute place finie $\mu$ de $L$.
\end{remarque}

\subsection{Formes anciennes et nouvelles en $v_0$.}
\label{oldetnews}
\subsubsection{Définition}\label{corpscoeffs}
On appelle {\em corps de définition} pour le triplet $(K,J,\rho)$
tout corps de nombres $L \subset \C$ comme 
en~\ref{modelecorpsdenombres} et~\ref{decomp}.
On appelle {\em bonne place} relativement à $(K,J,\rho)$
toute place finie de $L$ ou $L_0$ ne divisant aucun nombre premier des
ensembles $S_i$, $i=1,\ldots,4$, définis en~\ref{defs},~\ref{modelegroupe},
\ref{stabilitereseau} et~\ref{produitscalaire}. 
\subsubsection{Notations.}
Dans ce paragraphe, on conserve toutes les notations de~\ref{ntp},
on fixe un corps de définition $L$ pour $(K,J,\rho)$ 
et une bonne place $\muz$ de $L_0$.
On fixe aussi une place finie $v_0 \not \in \Sigma_\muz$,
un sous-groupe compact maximal spécial $K'_\vz$ 
de $G_\vz$ tel que les points de $X_\vz$  fixés par $K_\vz$ et
$K'_\vz$ soient voisins, et 
nous notons $B_\vz$ le sous-groupe d'Iwahori $K_\vz \cap K'_\vz$. 
Nous reprenons les notations $X=G_\vz/K_\vz$, $X'=G_\vz/K'_\vz$ et 
$A=G_\vz/B_\vz$ du paragraphe~\ref{prelimslocaux} pour la place 
$v_0$. 

Nous posons enfin $K':=\prod_{v \neq \vz} K_v \times K'_\vz$ et 
$B:=\prod_{v \neq \vz} K_v \times B_\vz$.
Puisque $\vz \not\in \Sigma$, 
la représentation $J$ de $K$ induit une représentation de $K'$ et $B$)
que l'on désigne par la même lettre 
et le corps de nombres $L$ (resp. la place $\muz$ de $L\cap \R$)
satisfait aussi à~\ref{modelecorpsdenombres} (resp.~\ref{modelegroupe},
\ref{stabilitereseau} et~\ref{produitscalaire}) 
pour les triplets $(K',J,\rho)$ et $(B,J,\rho)$,
de sorte que $S_{K',J,\rho,L}$ et $S_{B,J,\rho,L}$  
(resp. $S_{K',J,\rho,\anneau_\mu}$ et $S_{B,J,\rho,\anneau_\mu}$)
sont des modèles de $S_{K',J,\rho,\C}$ et $S_{B,J,\rho,\C}$  
(resp. $S_{K',J,\rho,L_\mu}$ et $S_{B,J,\rho,L_\mu}$).

On désigne par $R$ une $\anneau_\muz$-algèbre. 

\subsubsection{Interprétation locale en termes d'arbre}
\label{interpretationlocale}
Soit $\Gamma$ l'image dans $G_\vz$ du sous groupe de $G(F)$ 
défini par l'intersection $G(F)\cap \prod_{v\neq v_0} K_v$ dans 
$G(\A_{F,f})$. 
Comme $G(F)G_\vz$ est dense dans $G(\A_{F,f})$, 
et puisque la restriction de la 
représentation $\rho_\muz$ à $G_{v_0}$ est triviale, 
on a l'égalité
\begin{eqnarray*} 
S_{K,J,\rho,R} &=& \Hom_{\Gamma} 
(J_L\otimes \rho_L,\CC(X,R)).
\end{eqnarray*}
Ceci nous permet voir $S_{K,J,\rho,R}$ comme le sous-module de 
$\CC(X,R) \otimes_{\anneau_\muz}(V_{\anneau_\muz}^\ast 
\otimes V_{\anneau_\muz}^\ast)$ des fonctions 
covariantes sous $\Gamma$. 
Lorsque $R$ est muni d'une involution prolongeant celle de $\anneau_\muz$, 
le produit hermitien sur $S_{K,J,\rho,R}$ 
est l'intégrale sur $\Gamma\backslash X$ du produit tensoriel
du produit 
$\CC(X,\Z) \times \CC(X,\Z) \to \CC(X,\Z)$ par le produit hermitien de 
$(V_{\anneau_\muz}^\ast 
\otimes W_{\anneau_\muz}^\ast)\otimes_{\anneau_\muz} R$.

De même, on a
\begin{eqnarray*}
S_{K',J,\rho,R} &=& \Hom_{\Gamma} 
(J_L\otimes \rho_L,\CC(X',R)),\\
S_{B,J,\rho,R} &=& \Hom_{\Gamma} (J_L\otimes \rho_L,\CC(A,R)).
\end{eqnarray*}

\subsubsection{Opérateurs entre espaces de formes automorphes}
\label{operateursglobaux}

Les $R$-morphismes $T$, $U_{1}$, $U_{2}$, 
$U$, $U$ et $T_{B}$ 
définis en~\ref{operateurs} 
entres les modules $\CC(X,R)$, $\CC(X',R)$ et $\CC(A,R)$
commutent à l'action de $G_{\vz}$, donc de $\Gamma$, 
sur $X$, $X'$, $A'$ respectivement. 
Les morphismes qu'ils induisent par tensorisation par 
$V_{\anneau_\muz}^\ast 
\otimes W_{\anneau_\muz}^\ast$
définissent donc des $R$-morphismes  
entre $S_{K,J,\rho,R}$, $S_{K',J,\rho,R}$ et $S_{B,J,\rho,R}$ que nous
notons par les mêmes lettres. 
Ainsi, l'opérateur 
$T : S_{K,J,\rho,R} \rightarrow S_{K,J,\rho,R}$ est simplement 
l'opérateur de Hecke $T_\vz \in \Hecke(G_\vz,K_\vz Z_\vz)$ 
défini en~\ref{operateurhecke} et~\ref{ntp}. 

Notons que $\Hecke^{\Sigma_\muz}$ 
(resp. $\Hecke^{\Sigma_\muz\cup \{v_0\}}$)
agit sur $S_{K,J,\rho,R}$ (resp. $S_{K',J,\rho,R}$ et $S_{B,J,\rho,R}$).
On fait de $S_{B,J,\rho,R}$ un $\Hecke^\Sigma_\muz$-module
en faisant agir $T$ par $T_B$. 

\begin{lemme}\label{lemmeog}
Les $R$-morphismes $T$, $T_B$, $U_1$, $U_2$, $U$, $U'$ satisfont 
les assertions du lemme~\ref{relations}. De plus, 
$U_1$ et $U_2$ (resp. $U$ et $U'$) sont des morphismes de 
$\Hecke^{\Sigma_\muz}$-modules 
(resp. $\Hecke^{\Sigma_\muz \cup\{v_0\}}$).
\end{lemme}
\begin{pf}
La première assertion du lemme résulte de~\ref{interpretationlocale}.
La deuxième résulte de ce que les éléments de $G_v$ et $G_\vz$
commutent dans $G(\A_{F,f})$ pour tout $v \neq v_0$. 
\end{pf}

\subsubsection{Formes anciennes et formes nouvelles}\label{oldetnew}

Nous définissons l'espace $O_{B,J,\rho,R}$ 
des {\it formes
anciennes} et l'espace $N_{B,J,\rho,R}$ 
des {\it formes nouvelles} à valeurs dans $R$.

On pose tout d'abord 
$$ O_{B,J,\rho,\anneau_\muz }:=  
(U_1\oplus U_2)(S_{K,J,\rho,\anneau_\muz}^2)^\sat 
\subset S_{B,J,\rho,\anneau_\muz},$$ 
$$N_{B,J,\rho,\anneau_\muz}:= O_{B,J,K,\anneau_\muz}^\bot \subset S_{B,J,\rho,\anneau_\muz}.$$

Puis, pour toute $\anneau_\muz$-algèbre $R$,
$$ O_{B,J,\rho,\anneau_\muz }:= O_{B,J,\rho,\anneau_\muz}
\otimes_{\anneau_\muz} R  \subset S_{B,J,\rho,R},$$ 
$$ N_{B,J,\rho,\anneau_\muz }:= N_{B,J,\rho,\anneau_\muz}
\otimes_{\anneau_\muz} R  \subset S_{B,J,\rho,R},$$ 

Si $R$ est une $L_\muz$-algèbre, alors on a l'égalité $ O_{B,J,\rho,R}:=  
(U_1\oplus U_2)(S_{K,J,\rho,R}^2)$ et, 
puisque le produit hermitien sur $S_{B,J,\rho,\C}$ est défini positif,
la décomposition
$S_{B,J,\rho,R}=O_{B,J,\rho,R} \oplus N_{B,J,\rho,R}$. 
Ces assertions ne sont plus vraies {\em a priori} si $R=\anneau_\muz$, 
par exemple. 

\subsection{Le théorème}
\begin{theoreme}\label{theoaugniv} 
Soit $(K,J,\rho)$ un niveau, un type, et un poids
comme dans~\ref{ntp}. Soit $L$ un corps de définition pour
$(K,J,\rho)$ et $\mu$ une bonne place de $L$ relativement à
$(K,J,\rho)$ au sens de la définition~\ref{corpscoeffs}.
Soit $\Sigma$ la réunion de 
l'ensemble des places où $K$ n'est pas hyperspécial et 
soit $\psi$ un caractère de $\Hecke^{\Sigma}$ tel que 
$S_{K,J,\rho,L}(\psi)$ soit non nul.
Soit enfin $\vz$ une place de $F$ inerte dans $E$, hors de $\Sigma$
et de caractéristique résiduelle distincte de celle de $\mu$, 
et $\lambda=\psi(T_\vz) \in \anneau_\mu$. 
Si $\displaystyle{\lambda \not \in \{q(q^3+1),-(q^3+1)\}}$, notons  
$c$ le plus petit
entier vérifiant 
$$c \geq \val_\mu(\lambda-q(q^3+1))/2,$$ 
et si $\lambda=-(q^3+1)$, soit $c=\val_\mu(q^3+1)$;
alors il existe une congruence modulo $\mu^c$ entre  
$O_{B,J,\rho,\anneau_\mu}(\psi)$ et $N_{B,J,\rho,\anneau_\mu}$
dans $S_{B,J,\rho,\anneau_\mu}$. 
\end{theoreme}
\subsubsection{Remarques}\label{remarques}
\begin{itemize}
\item[(1)] Si $\lambda=-(q^3+1)$, alors
$\val_\mu(q^3+1)^2 \geq \val_\mu(\lambda-q(q^3+1))$, 
de sorte que l'on peut toujours prendre pour $c$ la partie
entière supérieure de $\val_\mu(\lambda-q(q^3+1))/2$, ce qui donne
bien
le théorème~\ref{augniv} de l'introduction. Ce cas correspond aux
représentations endoscopiques non tempérées (cf.~\ref{nontemp}).
\item[(2)] 
Si $q^3+1$ n'est pas nul modulo $\mu$, i.e. si l'on est en
caractéristique normale, alors le théorème reste vrai en 
prenant $c=\val_\mu(\lambda-q(q^3+1))$
(voir~\ref{remarquefinpreuve}). 
 
\item[(3)]
Soit $\Hecke$ un anneau commutatif contenant $\Hecke^{\Sigma}$
et agissant sur $S_{K,J,\rho,L}$ et $S_{B,J,\rho,L}$ 
de manière compatible aux structures entières
et aux opérateurs $U_1$ et $U_2$ définis en~\ref{operateursglobaux}.
Dans l'énoncé du théorème, on peut remplacer le caractère $\psi$
par un caractère de $\Hecke$ tel que le 
sous-espace propre généralisé $S_{K,J,\rho,L}(\psi)$ 
(cf.~\ref{sepgeneralises}) est non nul.  On obtient ainsi des formes
nouvelles sur lesquelles on connaît l'action modulo $\mu^c$ de l'algèbre 
$\Hecke$ tout entière. Cette précision peut être utile dans
les applications, en particulier si on prend pour $\Hecke$ le produit
tensoriel de $\Hecke^\Sigma$, avec un ou plusieurs anneaux $\Zc_v$, $v
\in \Sigma$, où $\Zc_v$ est le centre de l'algèbre de Hecke 
du type $(K_v,J_v)$
(voir~\ref{remarquefinpreuve2}).

\item[(4)] On peut vraisemblablement étendre le théorème à toutes les
places $\mu$ de $L$ à condition de prendre pour $c$ le plus petit
entier 
plus grand que  $\val_\mu(\lambda-q(q^3+1))/2 - n_\mu$, où $n_\mu$ est
un entier dépendant de $\mu \in S$ (et de $(K,J,\rho)$) mais pas de
$v_0$. Les détails sont laissés au lecteur intéressé.
\end{itemize}


Le reste de cette partie est consacrée à la preuve du théorème précédent.
\subsubsection{Notations} On pose 
\begin{eqnarray*} m &:=& \val_{\mu} (\lambda+(q^3+1)),\\  
n &:=& \val_{\mu} (\lambda-q(q^3+1))\\
d &:=& \dim_L S_{K,J,\rho,L}(\psi) \end{eqnarray*}

On note $\muz$ la trace de $\mu$ sur $L_0=L\cap \R$ et
$\anneau_\muz$ le complété de $\anneau_L$ en $\muz$, 
et on pose 
$$ M_0:=(U_1 \oplus U_2)(S_{K,J,\rho,\anneau_\muz}(\psi)^2)$$
et 
$$ M:=(U_1\oplus
U_2)(S_{K,J,\rho,\anneau_\mu}(\psi)^2)
=M_0 \otimes_{\anneau_\muz} \anneau_\mu.$$
Par le lemme~\ref{lemmeog}, 
on a $M_0^\sat = O_{B,J,\rho,\anneau_\muz}(\psi)$ et
$M^\sat = O_{B,J,\rho,\anneau_\mu}(\psi)=M_0^\sat 
\otimes_{\anneau_\muz} \anneau_\mu$.

Comme en~(\ref{phc}), 
on note $M_0^\anti$ (resp. $(M_0^\sat)^\anti$) l'ensemble des formes  
$\anneau_\muz$-antilinéaires sur $M_0$ (resp. $M_0^\sat$)
et $p_{M_0}\colon M_0 \to M_0^\anti$ 
(resp. $p_{M_0^\sat}$) la restriction du produit hermitien
de $S_{B,J,\rho,\anneau_\muz}$ à $M_0$ (resp. $M_0^\sat$).

\begin{lemme}\label{U1U2}\

(i) Si $\lambda\neq q(q^3+1),-(q^3+1)$, alors le 
morphisme de $L$-espaces vectoriels 
$U_1 \oplus U_2 : S_{K,J,\rho,L}(\psi)^2 \rightarrow 
O_{B,J,\rho,L}$ est injectif. En d'autres termes, les $\anneau_\mu$-réseaux 
$M$ et $M^\sat$ de $O_{K,J,\rho,L_\mu}(\eta)$ sont de rang $2d$. 

(ii) Si  $\lambda=-(q^3+1)$, alors l'image de  
$U_1 \oplus U_2 : S_{K,J,\rho,L}(\psi)^2 \rightarrow 
O_{B,J,\rho,L}$ est celle du morphisme injectif $U_1$
et les deux réseaux $M$ et $M^\sat$ de rang $d$ sont égaux. 
\end{lemme}
\begin{pf} Soient $(f_1,f_2)\in 
\ker (U_1 \oplus U_2)\cap S_{K,J,\rho,L}(\psi)$ 
et $u \in V_L\otimes W_L$. 
Alors $f_1(u)$ et
$f_2(u)$ appartiennent à $\CC(X,L)$ et vérifient 
$U_1 f_1(u)+U_2 f_2(u)=0$ et $T f_i(u) = \lambda
f_i(u)$ pour $i=1,2$. 

Dans le cas (i), le
lemme~\ref{relations} (iv), appliqué à $R=L$, implique que $f_1(u)$ et
$f_2(u)$ sont nulles. 
Comme ceci vaut pour tout $u \in V_L \otimes W_L$, c'est que
$f_1$ et $f_2$ sont nulles d'où la première partie du lemme.

Dans le cas (ii), alors le
lemme~\ref{relations}, (i), implique l'égalité
$U'Uf_i(u)=0$ pour $i=1,2$. Par l'assertion (ii) du même lemme, et 
puisque le produit hermitien sur $S_{K',J,\rho,L}$ est défini positif, 
on a donc $Uf_1(u)=Uf_2(u)=0$. Enfin, l'assertion (iii)  
montre que $f_1(u)=f_2(u)$, ce qui prouve que $M$ est l'image de $U_1$. 
Comme le lemme~\ref{relations}, (v), montre que $U_1$ est
injectif sur $S_{K,J,\rho,R}$ pour toute $\anneau_\mu$-algèbre $R$, 
le lemme~\ref{noyau} donne enfin l'égalité $M=M^\sat$.  
\end{pf}
\begin{lemme}\label{unautrelemme}
Si $\lambda \neq q(q^3+1), -(q^3+1)$, alors 
$\lng_{\anneau_\mu} (M^\sat/M) \leq dm/2$.
\end{lemme}
\begin{pf}
Le lemme~\ref{noyau} implique que pour tout entier $\alpha$ assez grand, 
le $\anneau_\mu$-module $M^\sat/M$ est isomorphe au noyau de
$(U_1 \oplus U_2)_{|S_{K,J,\rho,R}(\psi)^2}$
avec $R=\anneau/\mu^\alpha$.
D'après le lemme~\ref{relations}, (iii), 
ce module est encore isomorphe à 
$$\{f \in S_{K,J,\rho,R}(\psi), Uf \text{ est constante}\}$$
Par le lemme~\ref{lemmeog}, 
on a pour tout $f  \in S_{K,J,\rho,R}(\psi)$ telle que $Uf$ est
constante et tout $\tau \in \Hecke^{\Sigma_{\muz}}$ l'égalité
$$ (\eta_\const(\tau)- \psi(\tau)) Uf =  \tau Uf -\psi(\tau)Uf = 
U \tau f -\psi(\tau) U f=0 .$$ 
Or, comme les caractères $\eta_\const$ et $\psi$ ne sont pas 
congrus modulo $\mu$, le lemme~\ref{egalitecarac}
implique que leurs restrictions à $\Hecke(G_v,K_v)$ sont distinctes
pour une infinité de places $v \not \in \Sigma$. 
L'ensemble $\Sigma_\muz$ étant fini, on a nécessairement $Uf=0$.
 
On a donc prouvé que le $\anneau_\mu$-module  
$M^\sat/M$ est isomorphe à $\ker U_{|S_{K,J,\rho,R}(\psi)}$.
Or, le lemme~\ref{relations}, (i) et (ii), implique l'identité 
$U^\ast U = T+(q^3+1) = \lambda +(q^3+1)$. 
Puisque $\val_\mu(\lambda+(q^3+1))=m$, on a $\val_\mu(\det U^\ast U)= md$, 
et le lemme découle du lemme~\ref{longueur}.
\end{pf}

\begin{lemme}\label{unlemme} 

(i) Si $\lambda\neq q(q^3+1),-(q^3+1)$, alors 
$$\lng_{\anneau_\mu}((M^\anti_0/ p_{M_0}(M_0)) 
\otimes_{\anneau_{\muz}} \anneau_\mu)= d(m+n).$$

(ii) Si $\lambda=-(q^3+1)$, alors 
$\lng_{\anneau_\mu}((M^\anti_0/ p_{M_0}(M_0)) 
\otimes_{\anneau_{\muz}} \anneau_\mu)= d\val_\mu(q^3+1)$.
\end{lemme}
\begin{pf} (i) Par le lemme~\ref{noyau}, $\lng_{\anneau_\mu}(M_0^\anti / 
p_{M_0}(M_0))$
est la $\mu$-valuation du déterminant de $p_{M_0}$. 
Soit $e_1,\dots,e_d$ est une base orthonormée de 
$S_{K,J,\rho,\anneau_\muz}(\psi)$. Par le lemme~\ref{U1U2}, la famille 
$U_1(e_1),U_2(e_1),\dots,U_1(e_d),U_2(e_d)$ de vecteurs de $M_0$
en est une base sur laquelle la matrice de $p_{M_0}$ 
est, d'après la proposition~\ref{ii}, la matrice $(2d,2d)$ diagonale en des blocs $(2,2)$ 
tous égaux à
$$ \left( \begin{array}{cc} q^3+1 & \lambda \\ 
\lambda & q(q^3+1) + (q-1) \lambda \end{array} \right). $$ 
Le point (i) résulte donc du fait que le déterminant 
de la matrice (2,2) ci-dessus a pour $\mu$-valuation $n+m$
(voir la remarque suivant la proposition~\ref{ii}).

Pour prouver (ii), remarquons que, par le lemme~\ref{U1U2} 
et la proposition~\ref{ii}, la matrice de $p_{M_0}$ est
$q^3+1$ fois la matrice identité de dimension $d$.
\end{pf}

\begin{prop}\label{existencecong} Le module 
$(M_0^\sat)^\anti /p_{M_0^\sat}(M_0^\sat) $ 
contient un sous-module isomorphe à $\anneau_\mu/\mu^{c}$.
\end{prop}
\begin{pf} Notant $i$ l'injection canonique de $M_0$ dans $M_0^\sat$, 
et $i^\ast$ son adjoint pour $p_{M_0}$ et $p_{M_0^\sat}$, 
on a la suite d'injections 
$$\xymatrix{M_0 \ar[r]^i &  M_0^\sat \ar[r]^{p_{M_0^\sat}}  
& (M_0^\sat)^\anti \ar[r]^{i^\ast}& M_0^\anti },$$
dont la composée est $p_{M_0}$. 
Si $\lambda\neq q(q^3+1)$, alors 
$(M_0^\sat/M_0)\otimes \anneau_\mu= M^\sat/M$, 
et les lemmes~\ref{unlemme}\,(i) et~\ref{unautrelemme} impliquent 
\begin{eqnarray*}
\lng_{\anneau_\mu}\Big((M_0^\sat/p_{M_0^\sat} (M_0^\sat))\otimes 
\anneau_\mu\Big) 
&=& \lng_{\anneau_\mu} ((M_0^\anti/p_{M_0}(M_0))\otimes \anneau_\mu) 
- 2 \,\lng_{\anneau_\mu} (M^\sat/M))\\
&\geq& d(n+m)-2(dm/2)\\
&=& dn.
\end{eqnarray*}
Comme  $(M_0^\sat)^\anti /p_{M_0^\sat} (M_0^\sat)$ 
est engendré par $2d$ éléments
comme $\anneau_\muz$-module, la proposition en découle.

Si $\lambda=-(q^3+1)$, alors le lemme~\ref{U1U2} montre que
$(M_0^\sat/M_0)\otimes \anneau_\mu=0$, et la proposition
découle donc du lemme~\ref{unlemme}\,(ii) et du fait
que $(M_0^\sat)^\anti /p_{M_0^\sat} (M_0^\sat)$ 
est engendré par $d$ éléments.
\end{pf}

\begin{lemme}\label{decompold}
On a la décomposition orthogonale~:
$$O_{B,J,\rho,\anneau_\muz} = 
\oplus_{\eta} O_{B,J,\rho,\anneau_\muz}(\eta) $$
\end{lemme}
\begin{pf}
L'inclusion $O_{B,J,\rho,\anneau_\muz} \supset
\oplus_{\eta} O_{B,J,\rho,\anneau_\muz}(\eta)$
et l'orthogonalité proviennent de la décomposition 
du paragraphe~\ref{decomp} et du lemme~\ref{lemmeog}.
Le lemme~\ref{defs} montre de plus que 
les caractères $\eta$ qui interviennent sont deux-à-deux distincts
modulo toutes les places divisant $\muz$. Le 
lemme de Nakayama permet alors d'en déduire l'autre inclusion. 
\end{pf}

\subsubsection{Preuve du théorème~\ref{theoaugniv}}\label{preuve}
Posons $A=S_{B,J,\rho,\anneau_\muz}(\psi)=M_0^\sat$ 
et $B=N_{B,J,\rho,\anneau_\muz}$. 
Par le lemme~\ref{decompold}, on a 
$$(A\oplus B)^\perp
\oplus  A = O_{B,J,\rho,\anneau_\muz}.$$ 
On peut donc appliquer le lemme~\ref{modprod} aux 
sous-$\anneau_\muz$-modules 
$A$ et $B$ de $S_{B,J,\rho,\anneau_\muz}$, dont
on déduit que 
$(A \oplus B)^\sat /A\oplus B$ 
est isomorphe
à $A^\anti/p_{A}(A)$. 
Or, par la proposition~\ref{existencecong}.
$A^\anti/p_{A}(A) \otimes \anneau_\muz$ 
contient un sous-module isomorphe
à $\anneau_\mu/\mu^c$.
Le lemme~\ref{modcong} appliqué aux sous-$\anneau_\mu$-modules 
$A \otimes \anneau_\mu$ et $B\otimes\anneau_\mu$ de
$S_{B,J,\rho,\anneau_\mu}$ achève la preuve. 

\subsubsection{Preuve de la remarque~\ref{remarques}\,(2)}\label{remarquefinpreuve} 
Si $q^3+1$ n'est pas nul modulo $\mu$, prouvons qu'on peut prendre
$c=n$. On peut supposer $n \geq 1$, auquel cas 
$\lambda \equiv q(q^3+1) \not \equiv -(q^3+1) \pmod{\mu}$. 
La formule~(\ref{eq}) 
dans la preuve duu lemme~\ref{relations}\,(iv)
montre alors que le noyau de l'opérateur 
$U_1 \oplus U_2$ avec $R=\anneau_\mu/\mu$ est constitué 
de fonctions constantes,
et est donc nul puisque $\eta_\const \not \equiv \psi$. 
Par le lemme~\ref{unlemme}, on en déduit
l'égalité $M=(U_1 \oplus
U_2)(S_{K,J,\rho,\anneau_\mu}(\psi)^2)=M^\sat$, c'est-à-dire $M_0^\sat \otimes
\anneau_\mu = M_0 \otimes \anneau_\mu$. 

Par ailleurs, la preuve du lemme~\ref{unlemme}\,(i) montre que 
la matrice du produit hermitien $p_{M_0}$ sur $M_0$ 
dans une base adéquate est diagonale par blocs (2,2)
tous égaux à $$ \left( \begin{array}{cc} q^3+1 & \lambda \\ 
\lambda & q(q^3+1) + (q-1) \lambda \end{array} \right). $$ 
Comme cette matrice (2,2) a un déterminant de $\mu$-valuation $m+n=n$ 
et un premier mineur inversible dans $\anneau_\mu$, 
ses diviseur élémentaires sur $\anneau_\mu$ 
sont les idéaux $\anneau_\mu$ 
et $\mu^n$. On en déduit  
que $(M_0)^\anti/p_{M_0}(M_0)$ contient
un sous-module isomorphe à $\anneau_\mu/\mu^n$, et qu'il en va de
même de $(M_0^\sat)^\anti/p_{M_0^\sat}(M_0^\sat)$. On conclut comme
en~\ref{preuve}.

\subsubsection{Preuve de la remarque~\ref{remarques}\,(3)}
\label{remarquefinpreuve2} 
Dans ce cas, la 
décomposition de $S_{K,J,\rho,L}$ 
en sous-espaces propres généralisées pour $\Hecke$ 
donnée par~\ref{sepgeneralises}
est plus fine que celle de~\ref{decomp} (notons que cette fois,
la notion de sous-espaces propres généralisés est nécessaire, 
car l'action de $\Hecke$ n'est pas semi-simple 
{\em a priori}) et les
énoncés~\ref{U1U2} à~\ref{existencecong} restent donc valables. 
Suivant le lemme~\ref{defs}, si 
l'on définit $S$ de sorte que pour 
toute place $\mu$ ne divisant aucun élément
de $S$, les caractères $\eta$
de $\Hecke$ intervenant dans cette décomposition sont distincts modulo $\mu$,  
alors la preuve du lemme~\ref{decompold} fonctionne sans changement.  

\section{Représentations automorphes et représentations galoisiennes}
\label{sectionrepresentations}

La représentation unitaire $L^2(G(F) \backslash G(\A_F))$ de
$G(\A_F)$ est somme directe orthogonale de représentations
unitaires irréductibles
$$L^2(G(F) \backslash G(\A_F)) \simeq \oplus m(\pi) \pi$$
et les représentations $\pi$ qui apparaissent dans 
cette somme avec une multiplicité
$m(\pi) \geq 1$ sont appelés les {\it
représentations automorphes} de $G$.
À cette décomposition correspond une décomposition (cf.~\ref{defB})
$$\BB = \oplus  m(\pi) \pi_{lisse},$$
où $\pi_{lisse}$ désigne l'ensemble (dense) des vecteurs lisses de $\pi$.

\subsection{Preuve du corollaire~\ref{coraugniv}}
 
\subsubsection{Des représentations aux formes propres}
\label{repform}

Soient $K$, $(J,V)$, $(\rho,W)$ comme en~\ref{ntp}
et soit $\pi$ une représentation automorphe telle que $\pi_\infty \simeq
\rho$ et $\Hom_{K}(J,\pi)\neq 0$. Alors $\Hom_{K\times G(\A_{F,\infty})}
(J \otimes \rho,\pi)= S_{K,J,\rho,\C}$ 
(en effet, pour tout $f \in \Hom_{K\times G(\A_{F,\infty})}
(J \otimes \rho,\pi)$, $(v,w) \in V \times W$, 
$f(v \otimes w)$ appartient à $\pi_{lisse}$, donc à $\BB$).
De plus, si $\eta$ désigne le caractère donnant l'action 
de $\Hecke^\Sigma$ sur $\pi^{\prod_{v \not \in \Sigma} K_v}$, 
alors $\Hom_{K\times G(\A_{F,\infty})}
(J \otimes \rho,\pi) \subset   S_{K,J,\rho,\C}(\eta)$.

\subsubsection{Des formes nouvelles propres aux représentations} 
\label{formrep}

Soient de plus $v_0$ et $B$ comme en~\ref{oldetnews}, 
$\Sigma'$ un ensemble fini de places de $F$ contenant $\Sigma$
et $\eta$ un caractère de $\Hecke^{\Sigma'}$.
Si $N_{B,J,\rho,\C}(\eta)$ est non nul, alors on peut construire
une représentation automorphe $\pi'$ telle que $\pi'{}^\infty
\simeq \rho$ et $\Hom_{K}(J,\pi') \neq 0$, telle que 
l'action de $\Hecke^{\Sigma'}$ sur 
$\pi^{\prod_{v \not \in \Sigma'} K_v}$ soit donnée par $\eta$, 
et vérifiant de plus $\pi'{}^{B}\neq 0$ 
et $\pi'^{K}=0$. 

Soit en effet $0 \neq f \in N_{B,J,\rho,\C}(\eta)$,
soit $(v,w) \in V \times W$ tel que $f(v \otimes w)\neq 0$ 
et soit $\pi$ la représentation unitaire engendré par $f(v \otimes
w)$ dans $L^2(G(F) \backslash G(\A_F))$. Alors $\pi$ a une
décomposition finie orthogonale (avec éventuellement des multiplicités)
$\pi=\oplus_{i=1}^l \pi_i$ à laquelle correspond une décomposition
$f(v \otimes w)=\sum_{i=1}^l f_i$, avec $f_i \in \pi_i^B$ pour 
$i=1,\ldots,l$. Comme $f \in N_{B,J,\rho,\C}$, 
il existe un indice $j$, $1 \leq j \leq l$, 
tel que $f_j$ n'appartienne pas à l'image de $U_1\oplus U_2$ 
(cf.~\S\ref{combinatoire} et définition~\ref{oldetnew}). 
Comme $\pi_j$ est irréductible, la proposition~\ref{oldu} 
implique que $\pi_j^K=0$. Les autres assertions sont claires.

\subsubsection{Preuve du corollaire~\ref{coraugniv}}
Soient $\pi$, $v_0$ et $\mu$ vérifiant les hypothèses du 
corollaire~\ref{coraugniv}, $\pi$ étant telle que $\pi_\infty \simeq
\rho$ et $\Hom_{K}(J,\pi)\neq 0$. Si $\eta$ désigne le caractère 
de $\Hecke^\Sigma$ attaché à $\pi$, alors $S_{K,J,\rho,\C}(\eta) \neq 0$ 
par~\ref{repform}. 
Si l'on suppose que $\eta(T_{v_0})\equiv q(q^3+1) \pmod{\mu}$, 
alors le théorème~\ref{augniv}, et le lemme de Deligne-Serre,
donnent l'existence d'un caractère $\eta'$ de $\Hecke^{\Sigma_\mu}$ 
congru à $\eta$ modulo $\mu$ et d'une forme 
$0 \neq g \in N_{B,J,\rho,L}(\eta')$ propre pour $\eta'$.
La construction~\ref{formrep} permet de conclure à l'existence de $\pi'$. 
\smallskip

En complément du corollaire, on obtient~:
\begin{prop} Gardons les notations du corollaire~\ref{coraugniv}.
Alors $\pi'_{v_0}$ est isomorphe soit à $\pi^s$, soit à $\St$. 
Si $q^3+1 \not \equiv 0  \pmod{\mu}$, on peut assurer que
 $\pi'_{v_0} \simeq \St$.
\end{prop}
\begin{pf} La première assertion résulte immédiatement du corollaire
et de~\ref{invariants}. Pour la seconde, notons que la preuve du
corollaire montre qu'on construit $\pi'$ telle que $T_B$ agit sur
$(\pi'_{v_0})^{B_{v_0}}$ par $\lambda \equiv q(q^3+1) \pmod{\mu}$. 
Or, d'après le 
lemme~\ref{vptb}, $T_B$ agit sur $(\pi^s)^B$ par $-(q^3+1)$, et donc on
ne peut avoir $\pi'_{v_0}\simeq \pi^s$.
\end{pf}

\subsection{Un lemme de théorie des représentations}
Le but de ce paragraphe est de montrer le résultat suivant, 
que nous utilisons dans la preuve du
théorème~\ref{thexistence} (\S\ref{preuveexistence}).  

\subsubsection{Notations.} 
Soit $L\subset \C$ un corps de nombres stable sous la conjugaison 
complexe et $\mu$ une place de $L$.
Soit $\tau$ une représentation de $\Gal(\bar E/E)$ 
de dimension $d \leq 3$ sur $L_\mu$, 
telle qu'il existe un réseau sur $\anneau_\mu$ stable sous $\tau$. Si
$\Lambda$ désigne un tel réseau, 
on note $\tau_\Lambda$ 
la représentation sur $\anneau_\mu$ 
d'espace $\Lambda$ induite par $\tau$.
\begin{prop}\label{antiautodual}
Si $\tau$ est telle que $\tau^c\simeq \tau^\ast$, 
alors elle admet un réseau stable $\Lambda$
tel que $\tau_\Lambda^c\simeq \tau_\Lambda^\ast$
après une extension des scalaires ramifiée de degré $\leq 3$. 
\end{prop}
\noindent

\subsubsection{} La preuve de la proposition~\ref{antiautodual} utilise 
les deux lemmes suivants, pour lesquels nous introduisons
quelques notations.

Soient $\X$ l'immeuble de Bruhat-Tits attaché à $\Gl_3(L_\mu)$, $X$
l'ensemble de ses sommets, $\Sc$ la partie de $\X$ 
fixe par $\tau$ et $S=\Sc \cap X$. 
Rappelons que les points de $X$ correspondent aux
classes d'homothéties $[\Lambda]$ de $\anneau_\mu$-réseaux $\Lambda$
dans $L_\mu^3$, et que ceux de $S$ sont les classes $[\Lambda]$ avec $\Lambda$
stable sous $\tau$ (\cite[lemme 3.1.2]{ribellaiche}).

\begin{lemme}\label{premierlemme} 
Il existe un automorphisme de complexe polysimplicial
$b$ de $\X$, qui laisse stable $\Sc$, et tel que, pour
 tout $[\Lambda] \in S$, l'égalité $b([\Lambda])=[\Lambda']$ implique 
$$\tau_\Lambda^c \simeq \tau_{\Lambda'}^\ast.$$
\end{lemme}
\begin{pf} Voir \cite[7.3.3]{joeletgaetan}.
\end{pf}
\begin{lemme} Si $\tau$ est irréductible et vérifie 
$\tau^c \simeq \tau^\ast$, alors il existe un point de $S$ fixe par $b$ 
après une extension des scalaires ramifiée de degré $\leq 3$. 
\end{lemme}
\begin{pf} Comme $\tau$ est irréductible, l'ensemble $S$ est fini 
d'après~\cite[prop. 3.2.1 et remarque suivante]{ribellaiche}. 
Par~\cite{bt}, il existe donc une facette $F \subset \Sc$ 
telle que $b(F)=F$. L'isobarycentre $x$ de $F$ vérifie $b(x)=x$ et il
est clair qu'après une extension ramifiée de
degré $1 + \dim F$ de $L_\mu$, on a $x \in S$. 
\end{pf}
 
\subsubsection{Preuve de la proposition~\ref{antiautodual}}
Si $\tau$ est irréductible, il suffit de combiner les
deux lemmes ci-dessus. Si $\tau \simeq \tau_1 \oplus \tau_2$, où 
$\tau_1$ et $\tau_2$ sont irréductibles de dimension $1$ et $2$ 
respectivement, alors on a $\tau_i^c \simeq \tau_i^\ast$ pour $i=1,2$
et on est ramené au cas où $\tau$ est irréductible. 
Enfin, si $\tau$ est somme de caractères, alors 
$\Sc=S$ est réduit à un seul point et le résultat résulte
directement du lemme~\ref{premierlemme}.

\subsection{Preuve du théorème~\ref{thexistence}}\label{preuveexistence}

\subsubsection{Représentations galoisiennes attachées aux représentations 
automorphes}

 \label{repattachee}

Soit $\pi$ une représentation automorphe pour $G$.
D'après les travaux de Rogawski, Blasius et Rogawski 
(cf. \cite[th. 1.9.1]{br}), on peut attacher à $\pi$ un corps de définition
 $L$, et un 
système compatible $\rho_\mu : \Gal(\bar E/E) \rightarrow \Gl(M_\mu)$ de 
représentations continues de $\Gal(\bar E/E)$
de dimension $3$ sur $L_\mu$, $\mu$ parcourant l'ensemble des places finies 
de $L$, qui vérifient :

{\it Pour toute place finie $\mu$ de $L$ et toute  place finie $v$ de $F$ non ramifiée 
dans $E$ où $\pi$ est non ramifiée,  et 
toute place $w$ de $E$ au-dessus de $v$,
$\rho_\mu$ est non ramifiée en $w$, et le polynôme caractéristique de
$\rho_\mu(\Frob_w)$ coïncide avec celui de $\mat_{\pi,w}$.}

Fixons dorénavant une place finie $\mu$ de $L$. 
Compte tenu de la forme des matrices de Hecke 
(\ref{heckeinerte} et \ref{heckedec}), il vient facilement
que la représentation $\rho_\mu$ vérifient 
$$\rho_\mu^c \simeq \rho_\mu^\ast.$$
Appliquant la proposition~\ref{antiautodual}, 
et en rempla\c cant éventuellement $L_\mu$ par une extension ramifiée de
degré $\leq 3$, on en déduit qu'il 
existe un réseau $\Lambda$ de $L_\mu^3$, stable 
sous $\rho_\mu$, tel que la représentation $\rho_\Lambda$ 
de dimension $3$ sur $\anneau_\mu$
d'espace $\Lambda$ induite par $\rho_\mu$ vérifie
$\rho_\Lambda^c\simeq \rho_\Lambda^\ast$.

\subsubsection{}(Le lecteur pourra comparer avec~\cite[partie 6]{clozel}.)
Choisissons une base de $\Lambda$, ce qui nous permet de voir $\rho_\Lambda$ 
comme un morphisme $\Gal(\bar E/E) \rightarrow \Gl_3(\anneau_\mu)$ et
rappelons que $\gamma$ est un élément d'ordre $2$ de 
$\Gal(\bar E/F)$ relevant la conjugaison de $\Gal(E/F)$, cf.~\ref{corps}.
Pour toute matrice $M$ de $\Gl_3(L_\mu)$, on pose $M^\ast := {}^t M^{-1}$.
L'existence d'un isomorphisme $\rho_\Lambda^c \simeq \rho_\Lambda^\ast$, 
se traduit par celle d'une matrice $A \in \Gl_3(\anneau_\mu)$ 
telle que pour tout $g \in \Gal(\bar E/E)$
\begin{eqnarray} \label{antiauto} 
\rho_\Lambda(\gamma g \gamma^{-1}) = A \rho_\Lambda(g)^\ast A^{-1}. 
\end{eqnarray}
Notons $H$ l'image de $\rho_\Lambda$, 
$C$ le groupe à deux éléments $\{1,c\}$ 
La relation~(\ref{antiauto}) permet de faire agir $C$ 
sur le sous-groupe $H$ de $\Gl_3(\anneau/\mu^n)$ 
en faisant opérer $c$ par $M \mapsto A M^\ast A^{-1}$. 
On définit $\tilde H:= H \sd C$ comme étant le produit semi-direct 
attaché à cette action. La relation~(\ref{antiauto}) implique  
que $\rho_\Lambda$ se prolonge en un morphisme 
$\tilde \rho_\Lambda \colon\Gal(\bar E/F) \rightarrow \tilde H$ 
en posant
$$\left\{ \begin{array}{rcl}
\tilde \rho_\Lambda(\sigma) & := & \rho_\Lambda(\sigma) 
\sd 1, \ \sigma \in \Gal(\bar E/E) \\
\tilde \rho_\Lambda(\gamma) & := & 1 \sd c 
\end{array}\right.$$

\subsubsection{} Soit $n \geq 1$ un entier. 
Notons $H_n$ la réduction de $H$ modulo $\mu^n$, 
c'est-à-dire l'image de $H$ par le morphisme 
$\Gl_3(\anneau_\mu) \rightarrow \Gl_3(\anneau_\mu/\mu^n)$, 
notons $\tilde H_n$ le produit semi-direct $H_n\sd C$
et $\tilde \rho_n\colon \tilde \Gal(\bar E/F) \to \tilde H_n$ 
le morphisme de groupes induit par $\tilde \rho_\Lambda$.  
Soit $l$ la caractéristique résiduelle de $\mu$ et 
$\omega_{l^n} : \Gal(\bar E/F) \rightarrow 
(\Z/l^n\Z)^\ast$ le caractère cyclotomique. Comme 
l'élément $(1 \sd c, -1)$ de $\tilde H_n \times \Z/l^n \Z$
est l'image de $\gamma$ par $\tilde \rho_n \times \omega_{l^n}$, 
le théorème de $\check{\text{C}}$ebotarev implique 
l'existence d'un ensemble
de densité analytique strictement positive de places $v$ de $F$ telles
que $\tilde \rho_n(\Frob_v)  = 1 \sd c$ et $\omega_{l^n}(\Frob_v)=-1$.
Or, toute place $v$ de $F$ décomposée dans $E$ vérifie 
$\Frob_v \in \Gal(\bar E/E)$ soit 
$\tilde \rho_n (\Frob_v) \in \Gl_3(\anneau_\mu/\mu^n) \sd 1$, 
d'où l'existence d'un ensemble 
de densité analytique strictement positive de 
places $v$ de $F$ inertes
dans $E$ telles que 
$\tilde \rho_n(\Frob_v)  = 1 \sd c$ et $\omega_{l^n}(\Frob_v)=-1$.

Soit $v$ une place de cet ensemble
et $q$ le cardinal résiduel de $F_v$; vérifions 
que $v$ vérifie la conclusion du théorème~\ref{thexistence}. 
D'une part, on a l'égalité $q+1\equiv 0 (\mod \mu^n)$ puisque $q
=\omega_{l^\infty}(\Frob_v) \equiv -1 \pmod{\mu^n}$. 
D'autre part, en notant $w$ la place de $E$ au dessus de $v$, 
on a $\Frob_w=\Frob_v^2$, d'où l'égalité 
dans $\Gl_3(\anneau_\mu/\mu^n)$
$$\rho_n(\Frob_w) = (\tilde \rho_n(\Frob_v))^2=(1\sd c)^2  = 1
= \left(\begin{array}{ccc} q^2 & & \\ & 1 & \\ & & q^{-2} 
\end{array} \right),$$
ce qui par~\ref{satakeinerte} implique que 
l'opérateur de Hecke $T_v$ agit sur
$\pi_v^{K_v}$ par $q(q^3+1) \pmod{\mu^n}$. 

\subsection{Applications aux représentations endoscopiques non tempérées}

\subsubsection{Rappels}
\label{nontemp}
Une représentation automorphe pour $G$ est dite {\it endoscopique non 
tempérée} si elle vérifie les propriétés équivalentes suivantes
\begin{itemize} 
\item[(a)] $\pi_v$ est isomorphe à $\pi^n$ en une place $v$ inerte dans $E$.
\item[(b)] $\pi_v$ est isomorphe à $\pi^n$ en presque toute les places
$v$ inertes dans $E$.
\item[(b')] $\pi_v$ est isomorphe à $\pi^n$ en presque toute les places
$v$ inertes où $\pi$ est non ramifiée.
\item[(c)] Il existe deux caractères de Hecke $\chi$ et $\phi$ de $E$,
vérifiant  $\chi^c = \chi^{-1}$ et $\phi^c = \phi^{-1}$,  
tels que pour toute place
$w$ de $E$ la matrice de Hecke $\mat_w(\pi)$ soit
$$\diag(\chi(\pig_w)\phi(\pig_w)|\pig_w|^{1/2},\phi(\pig_w),\chi(\pig_w)\phi(\pig_w)|\pig_w|^{-1/2}),$$
où $\pig_w$ désigne  une uniformisante de $E_w$.
\item[(d)] $\pi$ est non tempérée en presque toute place et est de
dimension infinie.
\end{itemize}

L'équivalence entre (a) et (b) résulte de~\cite[th. 13.3.6 (c)]{roglivre},
l'équivalence entre (b) et (c) de~\cite[pages 395-398]{rog3} ;
que (c) implique (d) est clair, car les matrices de Hecke des
représentations $1$-dimensionnelles sont de la forme 
$$\diag(\phi(\pig_w)|\pig_w|,\phi(\pig_w),\phi(\pig_w)|\pig_w|^{-1})$$
et (d) implique (c) résulte de la classification des $A$-paquets de
$G$ rappelée 
en~\cite[2.9]{Rog1}, de la construction de la
représentation galoisienne associée et des conjectures de Weil,
prouvées par Deligne.

On conjecture en général que les représentations automorphes
qui sont tempérées en presque toutes places le sont partout. Pour $G=\U(3)$,
nous montrerons ce résultat à l'aide du
théorème~\ref{augniv} dans un prochain travail. 
Sans recourir à ce résultat, convenons 
d'appeler abusivement {\it tempérées} 
les représentations automorphes qui le sont
en presque toute place.
 
\begin{theoreme} Soit $\pi$ une représentation automorphe non
tempérée et $L$ un corps de définition pour $\pi$. 
Pour presque toute place $\mu$ de $L$
il existe une représentation automorphe tempérée $\pi'$ telle que
$\pi' \equiv \pi \pmod{\mu}$. \label{thtemp}
\end{theoreme}

\begin{pf} 
Le corollaire~\ref{coraugniv} et le théorème~\ref{thexistence} 
 montrent pour presque toute place $\mu$ de $L$, 
l'existence d'une représentation $\pi'$, d'une place inerte
 $v_0$ où $\pi$ est non ramifiée telles que
$(\pi'_{v_0})^{B_{v_0}} \neq 0$ et $(\pi'_{v_0})^{K_{v_0}} =0$ vérifiant
$\pi'\equiv \pi \pmod{\mu}$.

Reste à montrer qu'on peut imposer que $\pi'$ soit
tempérée. 

La forme des matrices de Hecke interdit que $\pi'$ soit de dimension
1, quitte à exclure  encore un nombre fini de $l$. Supposons par
l'absurde que $\pi'$ est endoscopique non tempérée, et soit
$\chi'$,$\phi'$
les caractères de Hecke de $E$ attachée à $\pi'$. On a $\chi \equiv
\chi' \pmod{\mu}, \psi'\equiv \psi; \pmod{\mu}$. En particulier, les
signes des équations fonctionnelles $\epsilon(\chi)$ et
$\epsilon(\chi')$ sont les mêmes
(cf. \cite{deligne}). D'après~\cite[théorème 1.2]{rog3}, la
multiplicité de  $\pi$
(resp. $\pi'$) est donnée par
$$m(\pi)=1/2(1+\epsilon(\chi)(-1)^{N+\deg_\Q F})$$
(resp. $m(\pi')=1/2(1+\epsilon(\chi')(-1)^{N'+\deg_\Q F})$) où $N$ (resp. $N'$)
est le nombre de places inertes ou ramifiées où $\pi$ est du type
$\pi^s$ (cf. {\it loc. cit} pour le sens exact de ``de type $\pi^s$'').  
En la place $v_0$, on a $\pi'_{v_0} \simeq  \pi^s$ 
car $\pi'_{v_0} \simeq \St$ exclurait que $\pi'$ soit endoscopique
(cf.~\cite[th. 14.6.4]{roglivre}), tandis que $\pi_{v_0} \simeq \pi^n$. 
Au places $v \neq v_0$, un choix adéquat de poids $K$ et du type $J$
assure que $\pi'_v$  de type $\pi^s$ si et seulement si $\pi_v$
l'est. On a donc $N'=N+1$. Ceci exclut que $m(\pi) m(\pi') \geq 1$, et cette 
contradiction prouve le théorème.
\end{pf}

{\tiny{ }\  \\
Joël Bellaïche, Philippe Graftieaux\\
jbellaic@math.unice.fr, graftiea@math.unice.fr\\
Laboratoire J.A.Dieudonne \\ 
Universite de Nice Sophia-Antipolis\\
Parc Valrose\\
06108 Nice Cedex 2\\}

\end{document}